\documentclass[11pt]{article}
\usepackage{amsfonts}
\usepackage{graphics}
\usepackage{indentfirst}
\usepackage{color}
\usepackage{cite}
\usepackage{latexsym}
\usepackage[paper=a4paper, left=2.1cm, right=2.1cm, top=2.2cm, bottom=1.5cm, headheight=5.5pt, footskip=0.8cm, footnotesep=0.8cm, centering, includefoot]{geometry}
\usepackage{amsmath}
\allowdisplaybreaks
\usepackage{amssymb}
\usepackage[dvips]{epsfig}
\usepackage{amscd}

\newtheorem{theorem}{Theorem}[section]
\newtheorem{remark}{Remark}[section]
\newtheorem{definition}{Definition}[section]
\newtheorem{lemma}{Lemma}[section]

\DeclareMathOperator{\curl}{curl}
\DeclareMathOperator{\loc}{loc}

\DeclareMathOperator{\divv}{div}

\makeatletter
\@addtoreset{equation}{section}
\makeatother
\makeatletter
\@addtoreset{equation}{section}
\makeatother

\title{{On formation of singularity of the full compressible magnetohydrodynamic equations with zero heat conduction}
\thanks{Supported by Fundamental Research Funds for the Central Universities (No. XDJK2017C050),
China Postdoctoral Science Foundation (No. 2017M610579), and the Doctoral Fund of Southwest University (No. SWU116033).}
}

\author{Xin Zhong\thanks{School of Mathematics and Statistics, Southwest University, Chongqing 400715,
People's Republic of China ({\tt xzhong1014@amss.ac.cn}).
}
}
\date{ }

\begin{document}
\maketitle

\begin{abstract}
We are concerned with the formation of singularity and breakdown of strong solutions to the Cauchy problem of the three-dimensional full compressible magnetohydrodynamic equations with zero heat conduction. It is proved that for the initial density allowing vacuum, the strong solution exists globally if the deformation tensor $\mathfrak{D}(\mathbf{u})$ and the pressure $P$ satisfy $\|\mathfrak{D}(\mathbf{u})\|_{L^{1}(0,T;L^\infty)}
+\|P\|_{L^{\infty}(0,T;L^\infty)}<\infty$. In particular, the criterion is independent of the magnetic field. The logarithm-type estimate for the Lam{\'e} system and some delicate energy estimates play a crucial role in the proof.
\end{abstract}


Math Subject Classification: 76W05; 35B65

\section{Introduction}
Let $\Omega\subset\mathbb{R}^3$ be a domain, the motion of a viscous, compressible, and heat conducting magnetohydrodynamic (MHD) flow in $\Omega$ can be described by the full compressible MHD equations
\begin{align}\label{1.1}
\begin{cases}
\rho_{t}+\divv(\rho\mathbf{u})=0,\\
(\rho\mathbf{u})_{t}+\divv(\rho\mathbf{u}\otimes\mathbf{u})
-\mu\Delta\mathbf{u}
-(\lambda+\mu)\nabla\divv\mathbf{u}+\nabla P=
\mathbf{b}\cdot\nabla\mathbf{b}-\frac12\nabla|\mathbf{b}|^2,\\
c_{\nu}[(\rho\theta)_{t}+\divv(\rho\mathbf{u}\theta)]
+P\divv\mathbf{u}-\kappa\Delta\theta
=2\mu|\mathfrak{D}(\mathbf{u})|^2+\lambda(\divv\mathbf{u})^2+\nu|\curl\mathbf{b}|^2,\\
\mathbf{b}_{t}-\mathbf{b}\cdot\nabla\mathbf{u}
+\mathbf{u}\cdot\nabla\mathbf{b}+\mathbf{b}\divv\mathbf{u}=\nu\Delta\mathbf{b},\\
\divv \mathbf{b}=0.
\end{cases}
\end{align}
Here, $t\geq0$ is the time, $x\in\Omega$ is the spatial coordinate, and $\rho, \mathbf{u}, P=R\rho\theta\ (R>0), \theta, \mathbf{b}$ are the fluid density, velocity, pressure, absolute temperature, and the magnetic field respectively; $\mathfrak{D}(\mathbf{u})$ denotes the deformation tensor given by
\begin{equation*}
\mathfrak{D}(\mathbf{u})=\frac{1}{2}(\nabla\mathbf{u}+(\nabla\mathbf{u})^{tr}).
\end{equation*}
The constant viscosity coefficients $\mu$ and $\lambda$ satisfy the physical restrictions
\begin{equation}\label{1.2}
\mu>0,\ 2\mu+3\lambda\geq0.
\end{equation}
Positive constants $c_\nu,\ \kappa$, and $\nu$ are respectively the heat capacity, the ratio of the heat conductivity coefficient over the heat capacity, and the magnetic diffusive coefficient.

There is huge literature on the studies about the theory of
well-posedness of solutions to the Cauchy problem and the initial boundary value problem (IBVP) for the compressible MHD system due to the physical importance, complexity, rich phenomena and mathematical challenges, refer to \cite{HW2008,HW2010,LH2015,LSX2016,LXZ2013,FY2009,XH2017} and references therein.
However, many physical important and mathematical fundamental problems are still open due to the lack of smoothing mechanism and the strong nonlinearity. When the initial density allows vacuum, the local large strong solutions to Cauchy problem and initial-boundary-value problem of 3D full MHD equations have been obtained, respectively, by Fan-Yu \cite{FY2009} and Xi-Hao \cite{XH2017}.
For the global well-posedness of strong solutions, Li-Xu-Zhang \cite{LXZ2013} and L{\"u}-Shi-Xu \cite{LSX2016} established the global existence and uniqueness of strong solutions to the 3D and 2D MHD equations, respectively, provided the smooth initial data are of small total energy. Furthermore, Hu-Wang \cite{HW2008,HW2010} showed the global existence of renormalized solutions to the compressible MHD equations for general large initial data. Nevertheless, it is an outstanding challenging open problem to investigate the global well-posedness for general large strong solutions with vacuum.

Therefore, it is important to study the mechanism of blow-up and structure of possible singularities of strong (or classical) solutions to the compressible MHD equations. The pioneering work can be traced to \cite{HX2005}, where He and Xin proved Serrin's criterion for strong solutions to the incompressible MHD system, that is,
\begin{equation}\label{1.3}
\lim_{T\rightarrow T^{*}}\|\mathbf{u}\|_{L^s(0,T;L^r)}=\infty,\ \text{for}\ \frac2s+\frac3r=1,\ 3<r\leq\infty,
\end{equation}
here $T^*$ is the finite blow up time. For the compressible isentropic MHD system, Xu-Zhang \cite{XZ2012} obtained the following Serrin type criterion
\begin{equation}\label{1.4}
\lim_{T\rightarrow T^{*}}\left(\|\divv\mathbf{u}\|_{L^1(0,T;L^\infty)}+\|\mathbf{u}\|_{L^s(0,T;L^r)}\right)=\infty,
\end{equation}
where $r$  and $s$ as in \eqref{1.3}. This criterion is similar to \cite{HLX20112} for 3D compressible isentropic Navier-Stokes equations, which shows that the mechanism of blow-up is independent of the magnetic field. For the 3D full compressible MHD system, Lu-Du-Yao \cite{LDY2011} proved that
\begin{equation}\label{1.5}
\lim_{T\rightarrow T^{*}}\left(\|\nabla\mathbf{u}\|_{L^1(0,T;L^\infty)}+\|\theta\|_{L^\infty(0,T;L^\infty)}\right)
=\infty
\end{equation}
under the assumption
\begin{equation}\label{1.6}
\mu>4\lambda.
\end{equation}
Recently, for the Cauchy problem and the IBVP of 3D full compressible MHD system, Huang-Li \cite{HL2013} proved that
\begin{equation}\label{1.7}
\lim_{T\rightarrow T^{*}}\left(\|\rho\|_{L^\infty(0,T;L^\infty)}+\|\mathbf{u}\|_{L^s(0,T;L^r)}\right)
=\infty,\ \text{for}\ \frac2s+\frac3r\leq1,\ 3<r\leq\infty.
\end{equation}
For more information on the blow-up criteria of compressible flows, we refer to \cite{HL2013,HLW2013,HLX20112,SWZ20112,WZ2014,X1998,XY2013, XZ2012,Z2017,DW2015,W2016,LCH2016,FLN2015,Z20172} and the references therein.

It should be noted that all the results mentioned above on the blow-up of strong (or classical) solutions of viscous, compressible, and heat conducting MHD flows are for $\kappa>0$. Very recently, Huang-Xin \cite{HX2016} obtained blow-up criteria for the non-isentropic compressible Navier-Stokes equations without heat-conductivity.
Therefore, it seems to be an interesting question to ask what the blow-up criterion is for the the system \eqref{1.1} with zero heat conduction.
In fact, this is the main aim of this paper.

When $\kappa=0$, and without loss of generality, take $c_\nu=R=1$, the system \eqref{1.1} can be written as
\begin{align}\label{1.10}
\begin{cases}
\rho_{t}+\divv(\rho\mathbf{u})=0,\\
(\rho\mathbf{u})_{t}+\divv(\rho\mathbf{u}\otimes\mathbf{u})
-\mu\Delta\mathbf{u}
-(\lambda+\mu)\nabla\divv\mathbf{u}+\nabla P=
\mathbf{b}\cdot\nabla\mathbf{b}-\frac12\nabla|\mathbf{b}|^2,\\
P_{t}+\divv(P\mathbf{u})
+P\divv\mathbf{u}
=2\mu|\mathfrak{D}(\mathbf{u})|^2+\lambda(\divv\mathbf{u})^2+\nu|\curl\mathbf{b}|^2,\\
\mathbf{b}_{t}-\mathbf{b}\cdot\nabla\mathbf{u}
+\mathbf{u}\cdot\nabla\mathbf{b}+\mathbf{b}\divv\mathbf{u}=\nu\Delta\mathbf{b},\\
\divv \mathbf{b}=0.
\end{cases}
\end{align}
The present paper is aimed at giving a blow-up criterion of strong solutions to the Cauchy problem  of the system \eqref{1.10} with the initial condition
\begin{equation}\label{1.11}
(\rho,\mathbf{u},P,\mathbf{b})(x,0)=(\rho_0,\mathbf{u}_0,P_0,\mathbf{b}_0)(x),\ \ x\in\mathbb{R}^3,
\end{equation}
and the far field behavior
\begin{equation}\label{1.12}
(\rho,\mathbf{u},P,\mathbf{b})(x,t)\rightarrow(0, \mathbf{0},0, \mathbf{0}),\ \text{as}\ |x|\rightarrow+\infty,\ t>0.
\end{equation}

Before stating our main result, we first explain the notations and conventions used throughout this paper. We denote by
\begin{equation*}
\int\cdot\text{d}x=\int_{\mathbb{R}^3}\cdot\text{d}x.
\end{equation*}
For $1\leq p\leq\infty$ and integer $k\geq0$, the standard Sobolev spaces are denoted by:
\begin{equation*}
\begin{split}
\begin{cases}
L^p=L^p(\mathbb{R}^3),\ W^{k,p}=W^{k,p}(\mathbb{R}^3), \ H^{k}=H^{k,2}(\mathbb{R}^3), \\
D_{0}^{1}=\{u\in L^6|\nabla u\in L^2\},\
D^{k,p}=\{u\in L_{\loc}^1|\nabla^k u\in L^p\}.
\end{cases}
\end{split}
\end{equation*}
Now we define precisely what we mean by strong solutions to the problem \eqref{1.10}--\eqref{1.12}.
\begin{definition}[Strong solutions]\label{def1}
$(\rho,\mathbf{u},P,\mathbf{b})$ is called a strong solution to \eqref{1.10}--\eqref{1.12} in $\mathbb{R}^3\times(0,T)$, if for some $q_0>3$,
\begin{equation*}
\begin{split}
\begin{cases}
\rho\geq0,\ \rho\in C([0,T];L^1\cap H^1\cap W^{1,q_0}),\ \rho_t\in C([0,T];L^{q_0}),\\
(\mathbf{u},\mathbf{b})\in C([0,T];D_{0}^{1}\cap D^{2,2})\cap L^{2}(0,T;D^{2,q_0}),\
\mathbf{b}\in C([0,T];H^2), \\
(\mathbf{u}_t,\mathbf{b}_t)\in L^{2}(0,T;D^{1,2}),\
(\sqrt{\rho}\mathbf{u}_{t},\mathbf{b}_t)\in L^{\infty}(0,T;L^{2}), \\
P\geq0,\ P\in C([0,T];L^1\cap H^1\cap W^{1,q_0}),\ P_t\in C([0,T];L^{q_0}), \\
\end{cases}
\end{split}
\end{equation*}
and $(\rho,\mathbf{u},P,\mathbf{b})$ satisfies both \eqref{1.10} almost everywhere in $\mathbb{R}^3\times(0,T)$ and \eqref{1.11} almost everywhere in $\mathbb{R}^3$.
\end{definition}

Our main result reads as follows:
\begin{theorem}\label{thm1.1}
 For constant $\tilde{q}\in(3,6]$, assume that the initial data $(\rho_0\geq0,\mathbf{u}_0,P_0\geq0,\mathbf{b}_0)$ satisfies
\begin{align}\label{A}
\begin{cases}
(\rho_0,P_0)\in L^1\cap H^1\cap W^{1,\tilde{q}},\ \mathbf{u}_0\in D_{0}^{1}\cap D^{2,2},\\
\sqrt{\rho_0}\mathbf{u}_0\in L^2,\ \mathbf{b}_0\in H^{2},\ \divv\mathbf{b}_0=0,
\end{cases}
\end{align}
and the compatibility conditions
\begin{equation}\label{A2}
-\mu\Delta\mathbf{u}_0-(\lambda+\mu)\nabla\divv\mathbf{u}_0+\nabla P_0-(\curl\mathbf{b}_0)\times\mathbf{b}_0=\sqrt{\rho_0}\mathbf{g}
\end{equation}
for some $\mathbf{g}\in L^2(\Omega)$.
Let $(\rho,\mathbf{u},P,\mathbf{b})$ be a strong solution to the problem \eqref{1.10}--\eqref{1.12}. If $T^{*}<\infty$ is the maximal time of existence for that solution,
then we have
\begin{align}\label{B}
\lim_{T\rightarrow T^{*}}\left(\|\mathfrak{D}(\mathbf{u})\|_{L^{1}(0,T;L^\infty)}
+\|P\|_{L^{\infty}(0,T;L^\infty)}\right)=\infty
\end{align}
provided that
\begin{equation}\label{C}
3\mu>\lambda.
\end{equation}
\end{theorem}

Several remarks are in order.
\begin{remark}\label{re1.1}
The local existence of a strong solution with initial data as in Theorem \ref{thm1.1} can be established in the same manner as \cite{FY2009,XH2017}. Hence, the maximal time $T^{*}$ is well-defined.
\end{remark}
\begin{remark}\label{re1.2}
It is worth noting that the blow-up criteria \eqref{B} is independent of the magnetic field. Moreover, compared with \cite{HX2016} for the non-isentropic Navier-Stokes equations with $\kappa=0$, according to \eqref{B}, the $L^\infty$ bound of the temperature $\theta$ is not the key point to make sure that the solution $(\rho,\mathbf{u},P,\mathbf{b})$ is a global one, and it may go to infinity in the vacuum region within the life span of our strong solution.
\end{remark}
\begin{remark}\label{re1.3}
In \cite{HX2016}, to obtain higher order derivatives of the solutions, the restriction $\mu>4\lambda$ plays a crucial role in the analysis. In fact, the condition $\mu>4\lambda$ is only used to get the upper bound of $\int\rho|\mathbf{u}|^rdx$ for some $r\geq4$ (see \cite[Lemma 3.2]{HX2016}).
Here, we derive the upper bound of $\int\rho|\mathbf{u}|^4dx$ under the assumption $3\mu>\lambda$ (see Lemma \ref{lem34}), which is weaker than $\mu>4\lambda$.
\end{remark}

We now make some comments on the analysis of this paper. We mainly make use of continuation argument to prove Theorem \ref{thm1.1}. That is, suppose that \eqref{B} were false, i.e.,
\begin{equation*}
\lim_{T\rightarrow T^*}\left(\|\mathfrak{D}(\mathbf{u})\|_{L^{1}(0,T;L^\infty)}
+\|P\|_{L^{\infty}(0,T;L^\infty)}\right)\leq M_0<\infty.
\end{equation*}
We want to show that
\begin{equation*}
\sup_{0\leq t\leq T^*}\left(\|(\rho,P)\|_{H^1\cap W^{1,\tilde{q}}}+\|\nabla\mathbf{u}\|_{H^1}
+\|\mathbf{b}\|_{H^2}\right) \leq C<+\infty.
\end{equation*}

Since the magnetic field is strongly coupled with the velocity field of the fluid in the compressible MHD system, some new difficulties arise in comparison with the problem for the compressible Navier-Stokes equations studied in \cite{HX2016}.
The following key observations help us to deal with the interaction of the magnetic field and the velocity field very well. First, we prove (see Lemma \ref{lem33}) that a control of $L_t^1L_x^\infty$-norm of the deformation tensor implies a control on the  $L_t^\infty L_x^\infty$-norm of the magnetic field $\mathbf{b}$. Then, motivated by \cite{WZ2014,HX2016}, we derive a priori estimates of the $L^2$-norm of $|\mathbf{u}||\nabla\mathbf{u}|$ in both space and time, which is the second key observation in this paper (see Lemma \ref{lem34}). Finally, the a priori estimates on the $L_t^\infty L_x^{\tilde{q}}$-norm of $(\nabla\rho,\nabla P)$ and the $L_t^1 L_x^{\infty}$-norm of the velocity gradient can be obtained (see Lemma \ref{lem37}) simultaneously by solving a logarithm Gronwall inequality based on a logarithm estimate for the Lam{\'e} system (see Lemma \ref{lem24}) and the a priori estimates we have derived.

The rest of this paper is organized as follows. In Section \ref{sec2}, we collect some elementary facts and inequalities that will be used later. Section \ref{sec3} is devoted to the proof of Theorem \ref{thm1.1}.

\section{Preliminaries}\label{sec2}

In this section, we will recall some known facts and elementary inequalities that will be used frequently later.

We begin with the following Gronwall's inequality, which plays a central role in proving a priori estimates on strong solutions $(\rho,\mathbf{u},P,\mathbf{b})$.
\begin{lemma}\label{lem21}
Suppose that $h$ and $r$ are integrable on $(a, b)$ and nonnegative a.e. in $(a, b)$. Further assume that $y\in C[a, b], y'\in L^1(a, b)$, and
\begin{equation*}
y'(t)\leq h(t)+r(t)y(t)\ \ \text{for}\ a.e\ t\in(a,b).
\end{equation*}
Then
\begin{equation*}
y(t)\leq \left[y(a)+\int_{a}^{t}h(s)\exp\left(-\int_{a}^{s}r(\tau)d\tau\right)ds\right]
\exp\left(\int_{a}^{t}r(s)ds\right),\ \ t\in[a,b].
\end{equation*}
\end{lemma}
{\it Proof.}
See  \cite[pp. 12--13]{T2006}.  \hfill $\Box$

Next, the following Gagliardo-Nirenberg inequality will be used later.
\begin{lemma}\label{lem22}
Let $1\leq p,q,r\leq\infty$, and $j,m$ are arbitrary integers satisfying $0\leq j<m$. Assume that $v\in C_{c}^\infty(\mathbb{R}^n)$. Then
\begin{equation*}
\|D^jv\|_{L^p}\leq C\|v\|_{L^q}^{1-a}\|D^mv\|_{L^r}^{a},
\end{equation*}
where
\begin{equation*}
-j+\frac{n}{p}=(1-a)\frac{n}{q}+a\left(-m+\frac{n}{r}\right),
\end{equation*}
and
\begin{equation*}
\begin{split}
a\in
\begin{cases}
[\frac{j}{m},1),\ \ \text{if}\ m-j-\frac{n}{r}\ \text{is an nonnegative integer},\\
[\frac{j}{m},1],\ \ \text{otherwise}.
\end{cases}
\end{split}
\end{equation*}
The constant $C$ depends only on $n,m,j,q,r,a$.
\end{lemma}
{\it Proof.}
See  \cite[Theorem]{N1959}.  \hfill $\Box$

Next, the following logarithm estimate will be used to estimate $\|\nabla\mathbf{u}\|_{L^\infty}$.
\begin{lemma}\label{lem23}
For $q\in(3,\infty)$, there is a constant $C(q)>0$ such that for all $\nabla\mathbf{v}\in L^2\cap D^{1,q}$, it holds that
\begin{equation}\label{2.2}
\|\nabla\mathbf{v}\|_{L^\infty}\leq C\left(\|\divv\mathbf{v}\|_{L^\infty}
+\|\curl\mathbf{v}\|_{L^\infty}\right)\log(e+\|\nabla^2\mathbf{v}\|_{L^q})
+C\|\nabla\mathbf{v}\|_{L^2}+C.
\end{equation}
\end{lemma}
{\it Proof.}
See  \cite[Lemma 2.3]{HLX20112}.  \hfill $\Box$

Finally, we consider the following Lam{\'e} system
\begin{align}\label{2.3}
\begin{cases}
-\mu\Delta\mathbf{U}-(\lambda+\mu)\nabla\divv\mathbf{U}=\mathbf{F},\ \ x\in\mathbb{R}^3,\\
\mathbf{U}\rightarrow\mathbf{0},\ \text{as}\ |x|\rightarrow\infty,
\end{cases}
\end{align}
where $\mathbf{U}=(U^1,U^2,U^3),\ \mathbf{F}=(F^1,F^2,F^3)$, and $\mu,\lambda$ satisfy \eqref{1.2}.

The following logarithm estimate for the Lam{\'e} system \eqref{2.3} will be used to estimate
$\|\nabla\mathbf{u}\|_{L^\infty}$ and
$\|\nabla\rho\|_{L^2\cap L^q}$.
\begin{lemma}\label{lem24}
Let $\mu,\lambda$ satisfy \eqref{1.2}. Assume that $\mathbf{F}=\divv\mathbf{g}$ where $\mathbf{g}=(g_{kj})_{3\times3}$ with
$g_{kj}\in L^2\cap L^r\cap D^{1,q}$ for $k,j=1,\cdots, 3, r\in(1,\infty)$, and $q\in(3,\infty)$. Then the Lam{\'e}
system \eqref{2.3} has a unique solution $\mathbf{U}\in D_0^1
\cap D^{1,r}\cap D^{2,q}$,
and there exists a generic positive constant $C$ depending only on $\mu, \lambda, q$, and $r$ such that
\begin{equation*}
\|\nabla\mathbf{U}\|_{L^r}\leq C\|\mathbf{g}\|_{L^r},
\end{equation*}
and
\begin{equation*}
\|\nabla\mathbf{U}\|_{L^\infty}
\leq C\left(1+\log(e+\|\nabla\mathbf{g}\|_{L^q})\|\mathbf{g}\|_{L^\infty}
+\|\mathbf{g}\|_{L^r}\right).
\end{equation*}
\end{lemma}
{\it Proof.}
See  \cite[Lemma 2.3]{HL2013}.  \hfill $\Box$

\section{Proof of Theorem \ref{thm1.1}}\label{sec3}

Let $(\rho,\mathbf{u},P,\mathbf{b})$ be a strong solution described in Theorem \ref{thm1.1}. Suppose that \eqref{B} were false, that is, there exists a constant $M_0>0$ such that
\begin{equation}\label{3.1}
\lim_{T\rightarrow T^*}\left(\|\mathfrak{D}(\mathbf{u})\|_{L^{1}(0,T;L^\infty)}
+\|P\|_{L^{\infty}(0,T;L^\infty)}\right)\leq M_0<\infty.
\end{equation}

First, the upper bound of the density could be deduced directly from \eqref{1.10}$_1$ and \eqref{3.1} (see \cite[Lemma 3.4]{HLX20112}).
\begin{lemma}\label{lem31}
Under the condition \eqref{3.1}, it holds that for any $T\in[0,T^*)$,
\begin{equation}\label{3.2}
\sup_{0\leq t\leq T}\|\rho\|_{L^1\cap L^\infty}\leq C,
\end{equation}
where and in what follows, $C,C_1,C_2$ stand for generic positive constants depending only on $M_0,\lambda,\mu,\nu,T^{*}$, and the initial data.
\end{lemma}

Next, we have the following standard estimate.
\begin{lemma}\label{lem32}
Under the condition \eqref{3.1}, it holds that for any $T\in[0,T^*)$,
\begin{equation}\label{3.3}
\sup_{0\leq t\leq T}\left(\|\sqrt{\rho}\mathbf{u}\|_{L^2}^2+\|P\|_{L^1\cap L^\infty}
+\|\mathbf{b}\|_{L^2}^2\right)+\int_{0}^T\left(\|\nabla\mathbf{u}\|_{L^2}^2
+\|\nabla\mathbf{b}\|_{L^2}^2\right)dt
\leq C.
\end{equation}
\end{lemma}
{\it Proof.}
It follows from \eqref{1.10}$_3$ that
\begin{equation}\label{7.01}
P_t+\mathbf{u}\cdot\nabla P+2P\divv\mathbf{u}
=F\triangleq2\mu|\mathfrak{D}(\mathbf{u})|^2
+\lambda(\divv\mathbf{u})^2+\nu|\curl\mathbf{b}|^2\geq0.
\end{equation}
Due to \eqref{3.1}, we can always define particle path before blowup time
\begin{align*}
\begin{cases}
\frac{d}{dt}\mathbf{X}(x,t)
=\mathbf{u}(\mathbf{X}(x,t),t),\\
\mathbf{X}(x,0)=x.
\end{cases}
\end{align*}
Thus, along particle path, we obtain from \eqref{7.01} that
\begin{align*}
\frac{d}{dt}P(\mathbf{X}(x,t),t)
=-2P\divv\mathbf{u}+F,
\end{align*}
which implies
\begin{equation*}
P(\mathbf{X}(x,t),t)=\exp\left(-2\int_{0}^t\divv\mathbf{u}ds\right)
\left[P_0+\int_{0}^t\exp\left(2\int_{0}^s\divv\mathbf{u}d\tau\right)Fds\right]\geq0.
\end{equation*}
As a result, we deduce from \eqref{3.1} that
\begin{equation}\label{3.03}
\sup_{0\leq t\leq T}\|P\|_{L^1\cap L^\infty}
\leq C.
\end{equation}

Multiplying \eqref{1.10}$_2$ and \eqref{1.10}$_4$ by $\mathbf{u}$ and $\mathbf{b}$ respectively, then adding the two resulting equations together, and integrating over $\mathbb{R}^3$, we obtain after integrating by parts that
\begin{align}\label{3.4}
& \frac{1}{2}\frac{d}{dt}\int\left(\rho|\mathbf{u}|^2
+|\mathbf{b}|^2\right)dx
+\int\left[\mu|\nabla\mathbf{u}|^2+(\lambda+\mu)(\divv\mathbf{u})^2
+\nu|\nabla\mathbf{b}|^2\right]dx \nonumber \\
& =\int P\divv\mathbf{u}dx \nonumber \\
& \leq(\lambda+\mu)\int (\divv\mathbf{u})^2dx+C(\lambda,\mu)\int P^2dx,
\end{align}
which combined with \eqref{3.03} gives
\begin{align}\label{3.04}
\frac{1}{2}\frac{d}{dt}\left(\|\sqrt{\rho}\mathbf{u}\|_{L^2}^2
+\|\mathbf{b}\|_{L^2}^2\right)
+\mu\|\nabla\mathbf{u}\|_{L^2}^2
+\nu\|\nabla\mathbf{b}\|_{L^2}^2\leq C
\end{align}
Integrating \eqref{3.04} with respect to $t$ and applying \eqref{3.03} lead to the desired \eqref{3.3}. This completes the proof of Lemma \ref{lem32}.
\hfill $\Box$

Inspired by \cite{HX2005},
we have the following the upper bound of the magnetic field $\mathbf{b}$.
\begin{lemma}\label{lem33}
Under the condition \eqref{3.1}, it holds that for any $T\in[0,T^*)$,
\begin{equation}\label{3.5}
\sup_{0\leq t\leq T}\|\mathbf{b}\|_{L^2\cap L^\infty}\leq C.
\end{equation}
\end{lemma}
{\it Proof.}
Multiplying \eqref{1.10}$_4$ by $q|\mathbf{b}|^{q-2}\mathbf{b}\ (q\geq2)$ and integrating the resulting equation over $\mathbb{R}^3$, we derive
\begin{align}\label{3.6}
\frac{d}{dt}\int|\mathbf{b}|^qdx
+\nu q(q-1)\int|\mathbf{b}|^{q-2}|\nabla\mathbf{b}|^2dx
& \leq q\int\left(\mathbf{b}\cdot\nabla\mathbf{u}-\mathbf{u}\cdot\nabla\mathbf{b}
-\mathbf{b}\divv\mathbf{u}\right)\cdot|\mathbf{b}|^{q-2}\mathbf{b}dx
\nonumber\\
& =q\int\left(\mathbf{b}\cdot\mathfrak{D}(\mathbf{u})-\mathbf{u}\cdot\nabla\mathbf{b}
-\mathbf{b}\divv\mathbf{u}\right)\cdot|\mathbf{b}|^{q-2}\mathbf{b}dx.
\end{align}
By the divergence theorem and \eqref{1.10}$_5$, we get
\begin{equation*}
-q\int(\mathbf{u}\cdot\nabla)\mathbf{b}\cdot|\mathbf{b}|^{q-2}\mathbf{b}dx
=\int\divv\mathbf{u}|\mathbf{b}|^qdx,
\end{equation*}
which together with \eqref{3.6} yields
\begin{align}\label{a}
\frac{d}{dt}\int|\mathbf{b}|^qdx
+\nu q(q-1)\int|\mathbf{b}|^{q-2}|\nabla\mathbf{b}|^2dx
\leq(2q+1)\|\mathfrak{D}(\mathbf{u})\|_{L^\infty}\int|\mathbf{b}|^qdx.
\end{align}
Consequently, from $q\geq2$ and \eqref{a}, we immediately have
\begin{align*}
\frac{d}{dt}\|\mathbf{b}\|_{L^q}
\leq\frac{2q+1}{q}\|\mathfrak{D}(\mathbf{u})\|_{L^\infty}\|\mathbf{b}\|_{L^q}
\leq3\|\mathfrak{D}(\mathbf{u})\|_{L^\infty}\|\mathbf{b}\|_{L^q}.
\end{align*}
Then Gronwall's inequality and \eqref{3.1} imply that for any $q\geq2$,
\begin{align}\label{a1}
\sup_{0\leq t\leq T}\|\mathbf{b}\|_{L^q}\leq C,
\end{align}
where $C$ is independent of $q$. Thus, letting $q\rightarrow\infty$ in \eqref{a1} leads to the desired \eqref{3.5} and finishes the proof of Lemma \ref{lem33}. \hfill $\Box$

Motivated by \cite{WZ2014,HX2016}, we can improve the basic estimate obtained in Lemma
\ref{lem31}.
\begin{lemma}\label{lem34}
Under the condition \eqref{3.1}, it holds that for any $T\in[0,T^*)$,
\begin{equation}\label{4.1}
\sup_{0\leq t\leq T}\|\rho^{1/4}\mathbf{u}\|_{L^4}^{4}
+\int_{0}^{T}\||\mathbf{u}||\nabla\mathbf{u}|\|_{L^2}^{2}
dt \leq C
\end{equation}
provided that
\begin{equation*}
3\mu>\lambda.
\end{equation*}
\end{lemma}
{\it Proof.}
Multiplying $\eqref{1.10}_2$ by $4|\mathbf{u}|^{2}\mathbf{u}$ and integrating the resulting equation over $\mathbb{R}^3$ yield that
\begin{align}\label{4.2}
 & \frac{d}{dt}\int \rho|\mathbf{u}|^4dx+4\int\left
[\mu|\mathbf{u}|^{2}|\nabla\mathbf{u}|^2+(\lambda+\mu)|\mathbf{u}|^{2}(\divv\mathbf{u})^2
+\frac{\mu|\nabla|\mathbf{u}|^2|^2}{2}\right]dx\nonumber\\
& =4\int\divv(|\mathbf{u}|^{2}\mathbf{u})Pdx
-8(\lambda+\mu)\int\divv\mathbf{u}|\mathbf{u}|\mathbf{u}\cdot\nabla|\mathbf{u}|dx \nonumber\\
& \quad +4\int|\mathbf{u}|^{2}\mathbf{u}\cdot\left(\divv(\mathbf{b}\otimes\mathbf{b})
-\nabla\left(\frac{|\mathbf{b}|^2}{2}\right)\right)dx \nonumber\\
& \leq 4\int\divv(|\mathbf{u}|^{2}\mathbf{u})Pdx
-8(\lambda+\mu)\int\divv\mathbf{u}|\mathbf{u}|\mathbf{u}\cdot\nabla|\mathbf{u}|dx \nonumber\\
& \quad +C\int|\mathbf{u}|^{2}|\nabla\mathbf{u}||\mathbf{b}|^2\mbox{d}x.
\end{align}
For the last term of the right-hand side of \eqref{4.2}, one obtains from H{\"o}lder's inequality, Sobolev's inequality, and \eqref{3.5} that, for any $\varepsilon_1\in (0,1)$,
\begin{align*}
C\int|\mathbf{u}|^{2}|\nabla\mathbf{u}||\mathbf{b}|^2dx
& \leq 4\mu\varepsilon_1\int|\mathbf{u}|^{2}|\nabla\mathbf{u}|^2dx
+C(\varepsilon_1)\int|\mathbf{u}|^{2}|\mathbf{b}|^4\mbox{d}x\\
& \leq 4\mu\varepsilon_1\int|\mathbf{u}|^{2}|\nabla\mathbf{u}|^2dx
+C(\varepsilon_1)\|\mathbf{u}\|_{L^6}^{2}\|\mathbf{b}\|_{L^{6}}^{4}\\
& \leq 4\mu\varepsilon_1\int|\mathbf{u}|^{2}|\nabla\mathbf{u}|^2dx
+C(\varepsilon_1)\|\nabla\mathbf{u}\|_{L^2}^{2},
\end{align*}
which together with \eqref{4.2} leads to
\begin{align}
& \frac{d}{dt}\int \rho|\mathbf{u}|^4\mbox{d}x+4\int\left[\mu(1-\varepsilon_1)|\mathbf{u}|^{2}|\nabla\mathbf{u}|^2
+(\lambda+\mu)|\mathbf{u}|^{2}(\divv\mathbf{u})^2
+\frac{\mu|\nabla|\mathbf{u}|^2|^2}{2}\right]dx  \nonumber\\
& \leq 4\int\divv(|\mathbf{u}|^{2}\mathbf{u})Pdx
-8(\lambda+\mu)\int\divv\mathbf{u}|\mathbf{u}|\mathbf{u}\cdot\nabla|\mathbf{u}|dx
+C(\varepsilon_1)\|\nabla\mathbf{u}\|_{L^2}^{2}.\nonumber
\end{align}
Consequently,
\begin{align}\label{4.3}
& \frac{d}{dt}\int \rho|\mathbf{u}|^4\mbox{d}x
+4\int_{\mathbb{R}^3\cap\{|\mathbf{u}|>0\}}
\left[\mu(1-\varepsilon_1)|\mathbf{u}|^{2}|\nabla\mathbf{u}|^2
+(\lambda+\mu)|\mathbf{u}|^{2}(\divv\mathbf{u})^2
+\frac{\mu|\nabla|\mathbf{u}|^2|^2}{2}\right]dx  \nonumber\\
& \leq 4\int_{\mathbb{R}^3\cap\{|\mathbf{u}|>0\}}\divv(|\mathbf{u}|^{2}\mathbf{u})Pdx
-8(\lambda+\mu)\int_{\mathbb{R}^3\cap\{|\mathbf{u}|>0\}}\divv\mathbf{u}|\mathbf{u}|\mathbf{u}
\cdot\nabla|\mathbf{u}|\mbox{d}x+C(\varepsilon_1)\|\nabla\mathbf{u}\|_{L^2}^{2}.
\end{align}
Direct calculations give that for $x\in \mathbb{R}^3\cap\{|\mathbf{u}|>0\},$
\begin{align}
&|\mathbf{u}|^2|\nabla\mathbf{u}|^2=|\mathbf{u}|^4
\left|\nabla\left(\frac{\mathbf{u}}{|\mathbf{u}|}\right)\right|^2
+|\mathbf{u}|^2|\nabla|\mathbf{u}||^2,\label{4.4}\\
&|\mathbf{u}|\divv \mathbf{u}=|\mathbf{u}|^2\divv\left(\frac{\mathbf{u}}{|\mathbf{u}|}\right)
+\mathbf{u}\cdot\nabla|\mathbf{u}|.\label{4.5}
\end{align}
For $\varepsilon_1,\varepsilon_2\in (0,1),$ we now define a nonnegative function as follows:
\begin{equation}\label{4.6}
k(\varepsilon_1,\varepsilon_2)=
\begin{cases}
\dfrac{\mu\varepsilon_2(3-\varepsilon_1)}{\lambda+\varepsilon_1\mu},
~~~\text{if}\ \ \lambda+\varepsilon_1\mu>0,\\
~~0,~~~~~~~~~~~~~~~\text{if}\ \ \lambda+\varepsilon_1\mu\leq0.
\end{cases}
\end{equation}

We prove \eqref{4.1} in two cases.

\textit{Case 1:} we assume that
\begin{equation}\label{4.7}
\int_{\mathbb{R}^3\cap\{|\mathbf{u}|>0\}}|\mathbf{u}|^4
\left|\nabla\left(\frac{\mathbf{u}}{|\mathbf{u}|}\right)\right|^2\mbox{d}x\leq k(\varepsilon_1,\varepsilon_2)\int_{\mathbb{R}^3\cap\{|\mathbf{u}|>0\}}
|\mathbf{u}|^2|\nabla|\mathbf{u}||^2dx.
\end{equation}
It follows from \eqref{4.3} that
\begin{equation}\label{4.8}
\frac{d}{dt}\int \rho|\mathbf{u}|^4dx
+4\int_{\mathbb{R}^3\cap\{|\mathbf{u}|>0\}}\Psi\mbox{d}x\leq 4\int_{\mathbb{R}^3\cap\{|\mathbf{u}|>0\}}\divv(|\mathbf{u}|^{2}\mathbf{u})P\mbox{d}x
+C(\varepsilon_1)\|\nabla\mathbf{u}\|_{L^2}^{2},
\end{equation}
where
\begin{align*}
\Psi=&\mu(1-\varepsilon_1)|\mathbf{u}|^2|\nabla\mathbf{u}|^2
+(\lambda+\mu)|\mathbf{u}|^2(\divv\mathbf{u})^2+2\mu|\mathbf{u}|^2|\nabla|\mathbf{u}||^2
+2(\lambda+\mu)\divv\mathbf{u}|\mathbf{u}|\mathbf{u}\cdot\nabla|\mathbf{u}|.
\end{align*}
Employing \eqref{4.4} and \eqref{4.5}, we find that
\begin{align*}
\Psi & =\mu(1-\varepsilon_1)|\mathbf{u}|^2|\nabla\mathbf{u}|^2
+(\lambda+\mu)|\mathbf{u}|^2(\divv\mathbf{u})^2+2\mu|\mathbf{u}|^2|\nabla|\mathbf{u}||^2\\
& \quad +2(\lambda+\mu)|\mathbf{u}|^2\divv\left(\frac{\mathbf{u}}{|\mathbf{u}|}\right)
\mathbf{u}\cdot\nabla|\mathbf{u}|+2(\lambda+\mu)|\mathbf{u}\cdot\nabla|\mathbf{u}||^2\\
&=\mu(1-\varepsilon_1)\left(|\mathbf{u}|^4
\left|\nabla\left(\frac{\mathbf{u}}{|\mathbf{u}|}\right)\right|^2
+|\mathbf{u}|^2|\nabla|\mathbf{u}||^2\right)+(\lambda+\mu)\left(|\mathbf{u}|^2\divv
\left(\frac{\mathbf{u}}{|\mathbf{u}|}\right)+\mathbf{u}\cdot\nabla|\mathbf{u}|\right)^2\\
& \quad +2\mu|\mathbf{u}|^2|\nabla|\mathbf{u}||^2+2(\lambda+\mu)|\mathbf{u}|^2
\divv\left(\frac{\mathbf{u}}{|\mathbf{u}|}\right)\mathbf{u}\cdot\nabla|\mathbf{u}|
+2(\lambda+\mu)|\mathbf{u}\cdot\nabla|\mathbf{u}||^2\\
&=\mu(1-\varepsilon_1)|\mathbf{u}|^4\left|\nabla
\left(\frac{\mathbf{u}}{|\mathbf{u}|}\right)\right|^2
+\mu(3-\varepsilon_1)|\mathbf{u}|^2|\nabla|\mathbf{u}||^2
-\frac{\lambda+\mu}{3}|\mathbf{u}|^4\left|\divv\left(\frac{\mathbf{u}}{|\mathbf{u}|}\right)
\right|^2\\
& \quad +3(\lambda+\mu)\left(\frac{2}{3}|\mathbf{u}|^2\divv\left(\frac{\mathbf{u}}{|\mathbf{u}|}
\right)+\mathbf{u}\cdot\nabla|\mathbf{u}|\right)^2\\
&\geq -(\lambda+\varepsilon_1\mu)|\mathbf{u}|^4
\left|\nabla\left(\frac{\mathbf{u}}{|\mathbf{u}|}\right)\right|^2+\mu(3-\varepsilon_1)|
\mathbf{u}|^2|\nabla|\mathbf{u}||^2.
\end{align*}
Here we have used the facts that $\lambda+\mu>0$\footnote{From \eqref{1.2} and $3\mu-\lambda>0$, we have $5\mu+2\lambda>0$. Then by \eqref{1.2} again one gets $7\mu+5\lambda>0$, which combined with \eqref{1.2} again implies $9\mu+8\lambda>0$. This together with \eqref{1.2} once more gives $11\mu+11\lambda>0$. Thus the result follows.} and
$$\left|\divv\left(\frac{\mathbf{u}}{|\mathbf{u}|}\right)\right|^2\leq 3\left|\nabla\left(\frac{\mathbf{u}}{|\mathbf{u}|}\right)\right|^2.$$
Then we derive from \eqref{4.7} and \eqref{4.6} that
\begin{align}\label{4.9}
4\int_{\mathbb{R}^3\cap\{|\mathbf{u}|>0\}}\Psi dx
&\geq [-4(\lambda+\varepsilon_1\mu)k(\varepsilon_1,\varepsilon_2)
+4\mu(3-\varepsilon_1)]
\int_{\mathbb{R}^3\cap\{|\mathbf{u}|>0\}}|\mathbf{u}|^2|\nabla|\mathbf{u}||^2dx\nonumber\\
&\geq 4\mu(3-\varepsilon_1)(1-\varepsilon_2)\int_{\mathbb{R}^3\cap\{|\mathbf{u}|>0\}}|\mathbf{u}|^2|\nabla|\mathbf{u}||^2dx.
\end{align}
Thus, substituting \eqref{4.9} into \eqref{4.8} and using \eqref{3.3}, \eqref{4.4}, and \eqref{4.7} yield
\begin{align*}
& \frac{d}{dt}\int \rho|\mathbf{u}|^4dx+4\mu(3-\varepsilon_1)(1-\varepsilon_2)
\int_{\mathbb{R}^3\cap\{|\mathbf{u}|>0\}}|\mathbf{u}|^2
|\nabla|\mathbf{u}||^2dx\\
&\leq 4\int_{\mathbb{R}^3\cap\{|\mathbf{u}|>0\}}\divv(|\mathbf{u}|^{2}\mathbf{u})Pdx+C(\varepsilon_1)
\|\nabla\mathbf{u}\|_{L^2}^{2}\\
&\leq C\int_{\mathbb{R}^3\cap\{|\mathbf{u}|>0\}}
|\mathbf{u}|^{2}|\nabla\mathbf{u}|Pdx+C(\varepsilon_1)
\|\nabla\mathbf{u}\|_{L^2}^{2}\\
&\leq \varepsilon\int_{\mathbb{R}^3\cap\{|\mathbf{u}|>0\}}
|\mathbf{u}|^{2}|\nabla\mathbf{u}|^2dx+C(\varepsilon)\|\mathbf{u}\|_{L^6}^2\|P\|_{L^3}^2
+C(\varepsilon_1)\|\nabla\mathbf{u}\|_{L^2}^{2}\\
&\leq\varepsilon(1+k(\varepsilon_1,\varepsilon_2))
\int_{\mathbb{R}^3\cap\{|\mathbf{u}|>0\}}|\mathbf{u}|^{2}
|\nabla|\mathbf{u}||^2dx+C(\varepsilon,\varepsilon_1)\|\nabla\mathbf{u}\|_{L^2}^{2}.
\end{align*}
Taking $\varepsilon=\dfrac{2\mu(3-\varepsilon_1)(1-\varepsilon_2)}{1+k(\varepsilon_1,\varepsilon_2)}$, we have
\begin{align*}
\frac{d}{dt}\int \rho|\mathbf{u}|^4\mbox{d}x
+2\mu(3-\varepsilon_1)(1-\varepsilon_2)
\int|\mathbf{u}|^2|\nabla|\mathbf{u}||^2\mbox{d}x
\leq C(\varepsilon_1,\varepsilon_2)\|\nabla\mathbf{u}\|_{L^2}^{2},
\end{align*}
which combined with \eqref{4.4} and \eqref{4.7} implies
\begin{align}\label{4.10}
\frac{\mbox{d}}{\mbox{d}t}\int \rho|\mathbf{u}|^4\mbox{d}x
+\varepsilon\int|\mathbf{u}|^2|\nabla\mathbf{u}|^2\mbox{d}x
\leq C(\varepsilon_1,\varepsilon_2)\|\nabla\mathbf{u}\|_{L^2}^{2}.
\end{align}

\textit{Case 2:} we assume that
\begin{equation}\label{4.11}
\int_{\mathbb{R}^3\cap\{|\mathbf{u}|>0\}}|\mathbf{u}|^4
\left|\nabla\left(\frac{\mathbf{u}}{|\mathbf{u}|}\right)\right|^2
dx > k(\varepsilon_1,\varepsilon_2)\int_{\mathbb{R}^3
\cap\{|\mathbf{u}|>0\}}|\mathbf{u}|^2|\nabla|\mathbf{u}||^2dx.
\end{equation}
It follows from \eqref{4.3} that
\begin{align*}
 & \frac{d}{dt}\int \rho|\mathbf{u}|^4\mbox{d}x+4\int_{\mathbb{R}^3\cap\{|\mathbf{u}|>0\}}
 \left[\mu(1-\varepsilon_1)|\mathbf{u}|^{2}
|\nabla\mathbf{u}|^2+(\lambda+\mu)|\mathbf{u}|^{2}(\divv\mathbf{u})^2
+\frac{\mu|\nabla|\mathbf{u}|^2|^2}{2}\right]dx
\nonumber\\
& \leq 4\int_{\mathbb{R}^3\cap\{|\mathbf{u}|>0\}}\divv(|\mathbf{u}|^{2}\mathbf{u})Pdx
-8(\lambda+\mu)\int_{\mathbb{R}^3\cap\{|\mathbf{u}|>0\}}\divv\mathbf{u}|\mathbf{u}|\mathbf{u}\cdot\nabla|\mathbf{u}|dx
+C(\varepsilon_1)\|\nabla\mathbf{u}\|_{L^2}^{2}\\
&\leq C\int_{\mathbb{R}^3\cap\{|\mathbf{u}|>0\}}P|\mathbf{u}|^{2}|\nabla\mathbf{u}|\mbox{d}x+4(\lambda+\mu)
\int_{\mathbb{R}^3\cap\{|\mathbf{u}|>0\}}|\mathbf{u}|^2|\nabla|\mathbf{u}||^2dx  \\
& \quad +4(\lambda+\mu)\int_{\mathbb{R}^3\cap\{|\mathbf{u}|>0\}}|\mathbf{u}|^{2}|
\divv\mathbf{u}|^2dx+C(\varepsilon_1)\|\nabla\mathbf{u}\|_{L^2}^{2},
\end{align*}
which implies that
\begin{align}\label{4.12}
& \frac{d}{dt}\int \rho|\mathbf{u}|^4dx+4\mu(1-\varepsilon_1)
\int_{\mathbb{R}^3\cap\{|\mathbf{u}|>0\}}|\mathbf{u}|^{2}
|\nabla\mathbf{u}|^2dx+4(\mu-\lambda)
\int_{\mathbb{R}^3\cap\{|\mathbf{u}|>0\}}|\mathbf{u}|^{2}
|\nabla|\mathbf{u}||^2dx   \nonumber\\
& \leq C\int_{\mathbb{R}^3\cap\{|\mathbf{u}|>0\}}P|\mathbf{u}|^{2}|\nabla\mathbf{u}|dx
+C(\varepsilon_1)\|\nabla\mathbf{u}\|_{L^2}^{2}.
\end{align}
Inserting \eqref{4.4} into \eqref{4.12} yields
\begin{align*}
& \frac{d}{dt}\int \rho|\mathbf{u}|^4dx+[8\mu-4(\varepsilon_1\mu+\lambda)]
\int_{\mathbb{R}^3\cap\{|\mathbf{u}|>0\}}|\mathbf{u}|^{2}
|\nabla|\mathbf{u}||^2\mbox{d}x \nonumber \\
&\quad +4\mu(1-\varepsilon_1)
\int_{\mathbb{R}^3\cap\{|\mathbf{u}|>0\}}|\mathbf{u}|^{4}
\left|\nabla\left(\frac{\mathbf{u}}{|\mathbf{u}|}\right)\right|^2dx\\
&\leq C\int_{\mathbb{R}^3\cap\{|\mathbf{u}|>0\}}P|\mathbf{u}|^{2}
|\nabla|\mathbf{u}||dx+C\int_{\mathbb{R}^3\cap\{|\mathbf{u}|>0\}}P|\mathbf{u}|^{3}
\left|\nabla\left(\frac{\mathbf{u}}{|\mathbf{u}|}\right)\right|\mbox{d}x
+C(\varepsilon_1)\|\nabla\mathbf{u}\|_{L^2}^{2}\\
& \leq C\int_{\mathbb{R}^3\cap\{|\mathbf{u}|>0\}}P|\mathbf{u}|^{2}|\nabla|\mathbf{u}||dx
+4\mu(1-\varepsilon_1)\varepsilon_3\int_{\mathbb{R}^3\cap\{|\mathbf{u}|>0\}}|\mathbf{u}|^{4}
\left|\nabla\left(\frac{\mathbf{u}}{|\mathbf{u}|}\right)\right|^2dx\\
& \quad +C(\varepsilon_1,\varepsilon_3)\|\mathbf{u}\|_{L^6}^2\|P\|_{L^3}^2
+C(\varepsilon_1)\|\nabla\mathbf{u}\|_{L^2}^{2}\\
&\leq C\int_{\mathbb{R}^3\cap\{|\mathbf{u}|>0\}}P|\mathbf{u}|^{2}|\nabla|\mathbf{u}||dx
+4\mu(1-\varepsilon_1)\varepsilon_3
\int_{\mathbb{R}^3\cap\{|\mathbf{u}|>0\}}|\mathbf{u}|^{4}
\left|\nabla\left(\frac{\mathbf{u}}{|\mathbf{u}|}\right)
\right|^2\mbox{d}x\\
& \quad +C(\varepsilon_1,\varepsilon_3)\|\nabla\mathbf{u}\|_{L^2}^{2}
\end{align*}
with $\varepsilon_3\in (0,1)$. Hence we have
\begin{align*}
& \frac{d}{dt}\int \rho|\mathbf{u}|^4dx+[8\mu-4(\lambda+\varepsilon_1\mu)]
\int_{\mathbb{R}^3\cap\{|\mathbf{u}|>0\}}|\mathbf{u}|^{2}|\nabla|\mathbf{u}||^2dx
\nonumber \\
&\quad+4\mu(1-\varepsilon_1)(1-\varepsilon_3)\int_{\mathbb{R}^3\cap\{|\mathbf{u}|>0\}}|\mathbf{u}|^{4}
\left|\nabla\left(\frac{\mathbf{u}}{|\mathbf{u}|}\right)\right|^2dx\\
&\leq C\int_{\mathbb{R}^3\cap\{|\mathbf{u}|>0\}}P|\mathbf{u}|^{2}
|\nabla|\mathbf{u}||\mbox{d}x+C(\varepsilon_1,\varepsilon_3)\|\nabla\mathbf{u}\|_{L^2}^{2}.
\end{align*}
This together with \eqref{4.11} and \eqref{3.1} leads to
\begin{align*}
& \frac{\mbox{d}}{\mbox{d}t}\int \rho|\mathbf{u}|^4\mbox{d}x+k_1(\varepsilon_1,\varepsilon_2,\varepsilon_3,\varepsilon_4)
\int_{\mathbb{R}^3\cap\{|\mathbf{u}|>0\}}|\mathbf{u}|^{2}|\nabla|\mathbf{u}||^2\mbox{d}x
\nonumber \\
&\quad +k_2(\varepsilon_1,\varepsilon_3,\varepsilon_4)
\int_{\mathbb{R}^3\cap\{|\mathbf{u}|>0\}}|\mathbf{u}|^{4}\left|\nabla\left(\frac{\mathbf{u}}{|\mathbf{u}|}\right)\right|^2\mbox{d}x\\
&\leq C\int_{\mathbb{R}^3\cap\{|\mathbf{u}|>0\}}P|\mathbf{u}|^{2}|\nabla|\mathbf{u}||\mbox{d}x
+C(\varepsilon_1,\varepsilon_3)\|\nabla\mathbf{u}\|_{L^2}^{2},
\end{align*}
where
\begin{align*}
& k_1(\varepsilon_1,\varepsilon_2,\varepsilon_3,\varepsilon_4)
=4\mu(1-\varepsilon_1)(1-\varepsilon_3)(1-\varepsilon_4)k(\varepsilon_1,\varepsilon_2)
+8\mu-4(\lambda+\varepsilon_1\mu),\\
& k_2(\varepsilon_1,\varepsilon_3,\varepsilon_4)
=4\mu(1-\varepsilon_1)(1-\varepsilon_3)\varepsilon_4,
\end{align*}
with $\varepsilon_i\in (0,1),~~i=1,2,3,4.$\\
Since $k_2(\varepsilon_1,\varepsilon_3,\varepsilon_4)>0$ for all $(\varepsilon_1,\varepsilon_3,\varepsilon_4)\in (0,1)\times(0,1)\times(0,1),$ we only need to show that there exists $(\varepsilon_1,\varepsilon_2,\varepsilon_3,\varepsilon_4)\in (0,1)\times(0,1)\times(0,1)\times(0,1)$ such that
$$k_1(\varepsilon_1,\varepsilon_2,\varepsilon_3,\varepsilon_4)>0.$$
In fact, if $\lambda<0,$ take $\varepsilon_1=-\frac{\lambda}{m\mu}\in(0,1)$, with the positive integer $m$ large enough, then we have $$\varepsilon_1\mu+\lambda=\frac{m-1}{m}\lambda<0,$$ which implies that $k(\varepsilon_1,\varepsilon_2)=0,$ and hence
$$k_1(\varepsilon_1,\varepsilon_2,\varepsilon_3,\varepsilon_4)
=8\mu-4(\lambda+\varepsilon_1\mu)>8\mu>0.$$
If $\lambda=0,$ then $\lambda+\varepsilon_1\mu>0,$ which implies that
$$k_1(\varepsilon_1,\varepsilon_2,\varepsilon_3,\varepsilon_4)
=\dfrac{4\mu(1-\varepsilon_1)(1-\varepsilon_3)(1-\varepsilon_4)
(3-\varepsilon_1)\varepsilon_2}{\varepsilon_1}+8\mu-4\varepsilon_1\mu> 4\mu>0.$$
If $0<\lambda<3\mu,$ then we have $\lambda+\varepsilon_1\mu>0$ and then
$$k_1(\varepsilon_1,\varepsilon_2,\varepsilon_3,\varepsilon_4)
=\dfrac{4\mu^2(1-\varepsilon_1)(1-\varepsilon_3)
(1-\varepsilon_4)(3-\varepsilon_1)\varepsilon_2}{\lambda+\varepsilon_1\mu}
+8\mu-4(\lambda+\varepsilon_1\mu).$$
Notice that $k_1(\varepsilon_1,\varepsilon_2,\varepsilon_3,\varepsilon_4)$ is continuous over $[0,1]\times[0,1]\times[0,1]\times[0,1],$ and
$$k_1(0,1,0,0)=\dfrac{12\mu^2}{\lambda}+8\mu-4\lambda
=4\lambda^{-1}(\lambda+\mu)(3\mu-\lambda)>0,$$
so there exists $(\varepsilon_1,\varepsilon_2,\varepsilon_3,\varepsilon_4)\in (0,1)\times(0,1)\times(0,1)\times(0,1)$ such that
$$k_1(\varepsilon_1,\varepsilon_2,\varepsilon_3,\varepsilon_4)>0.$$
So we have
\begin{align*}
&\frac{\mbox{d}}{\mbox{d}t}\int \rho|\mathbf{u}|^4\mbox{d}x+k_1(\varepsilon_1,\varepsilon_2,\varepsilon_3,\varepsilon_4)
\int_{\mathbb{R}^3\cap\{|\mathbf{u}|>0\}}|\mathbf{u}|^{2}|\nabla|\mathbf{u}||^2
dx \nonumber \\
& \quad +k_2(\varepsilon_1,\varepsilon_3,\varepsilon_4)\int_{\mathbb{R}^3\cap\{|\mathbf{u}|>0\}}
|\mathbf{u}|^{4}
\left|\nabla\left(\frac{\mathbf{u}}{|\mathbf{u}|}\right)\right|^2\mbox{d}x
\nonumber\\
&\leq \frac{k_1(\varepsilon_1,\varepsilon_2,\varepsilon_3,\varepsilon_4)}{2}
\int_{\mathbb{R}^3\cap\{|\mathbf{u}|>0\}}|\mathbf{u}|^{2}|\nabla|\mathbf{u}||^2dx \nonumber \\
& \quad +C(\varepsilon_1,\varepsilon_2,\varepsilon_3,\varepsilon_4)\|\mathbf{u}\|_{L^6}^2
\|P\|_{L^3}^2
+C(\varepsilon_1,\varepsilon_3)\|\nabla\mathbf{u}\|_{L^2}^{2}
\nonumber\\
&\leq \frac{k_1(\varepsilon_1,\varepsilon_2,\varepsilon_3,\varepsilon_4)}{2}
\int_{\mathbb{R}^3\cap\{|\mathbf{u}|>0\}}|\mathbf{u}|^{2}|\nabla|\mathbf{u}||^2dx +C(\varepsilon_1,\varepsilon_2,\varepsilon_3,\varepsilon_4)\|\nabla\mathbf{u}\|_{L^2}^{2}.
\end{align*}
Therefore,
\begin{align}\label{4.13}
 & \frac{d}{dt}\int \rho|\mathbf{u}|^4dx+\frac{k_1(\varepsilon_1,\varepsilon_2,\varepsilon_3,\varepsilon_4)}{2}
\int|\mathbf{u}|^{2}|\nabla|\mathbf{u}||^2dx+k_2(\varepsilon_1,\varepsilon_3,\varepsilon_4)
\int|\mathbf{u}|^{4}\left|\nabla\left(\frac{\mathbf{u}}{|\mathbf{u}|}\right)\right|^2dx
\nonumber\\
& \leq C(\varepsilon_1,\varepsilon_2,\varepsilon_3,\varepsilon_4)\|\nabla\mathbf{u}\|_{L^2}^{2}.
\end{align}
From \eqref{4.10}, \eqref{4.13}, and \eqref{4.4}, we conclude that if $3\mu>\lambda,$ there exists a constant $\bar{C}>0$ such that
\begin{align*}
\frac{d}{dt}\int \rho|\mathbf{u}|^4dx+\bar{C}\int|\mathbf{u}|^2|\nabla\mathbf{u}|^2dx
\leq C\|\nabla\mathbf{u}\|_{L^2}^{2},
\end{align*}
which together with \eqref{3.3} and Gronwall's inequality gives the desired \eqref{4.1}.
\hfill $\Box$

Let $E$ be the specific energy defined by
\begin{equation}\label{3.13}
E=\theta+\frac{|\mathbf{u}|^2}{2}.
\end{equation}
Let $G$ be the effective viscous flux, $\pmb\omega$ be vorticity given by
\begin{equation}\label{3.14}
G=(\lambda+2\mu)\divv\mathbf{u}-\left(P+\frac{|\mathbf{b}|^2}{2}\right),\
\pmb\omega=\curl\mathbf{u}.
\end{equation}
Then the momentum equations \eqref{1.10}$_2$ can be rewritten as
\begin{equation}\label{3.15}
\rho\dot{\mathbf{u}}-\mathbf{b}\cdot\nabla\mathbf{b}
=\nabla G-\curl\pmb\omega,
\end{equation}
where $\dot{\mathbf{u}}\triangleq\mathbf{u}_t+\mathbf{u}\cdot\nabla\mathbf{u}$.

The following lemma gives the estimates on the spatial gradients of both the velocity and the magnetic field, which are crucial for deriving the higher order estimates of the solution.
\begin{lemma}\label{lem35}
Under the condition \eqref{3.1}, it holds that for any $T\in[0,T^*)$,
\begin{equation}\label{5.1}
\sup_{0\leq t\leq T}\left(\|\nabla\mathbf{u}\|_{L^2}^{2}
+\|\nabla\mathbf{b}\|_{L^2}^{2}\right)
+\int_{0}^{T}\left(\|\sqrt{\rho}\dot{\mathbf{u}}\|_{L^2}^{2}
+\|\mathbf{b}_t\|_{L^2}^{2}+\|\nabla^2\mathbf{b}\|_{L^2}^2\right)dt \leq C.
\end{equation}
\end{lemma}
{\it Proof.}
Multiplying \eqref{1.10}$_2$ by $\mathbf{u}_{t}$ and integrating the resulting equation over $\mathbb{R}^3$ give rise to
\begin{align}\label{5.2}
&\frac{1}{2}\frac{d}{dt}\int\left(\mu|\nabla\mathbf{u}|^2
+(\lambda+\mu)(\divv\mathbf{u})^2
\right)dx+\int\rho|\dot{\mathbf{u}}|^2dx\nonumber\\
& = \int\rho\dot{\mathbf{u}}\cdot(\mathbf{u}\cdot\nabla)\mathbf{u}dx
+\int\left(P+\frac{|\mathbf{b}|^2}{2}\right)\divv\mathbf{u}_tdx
-\int(\mathbf{b}\otimes\mathbf{b}):\nabla\mathbf{u}_{t}dx
\nonumber\\
& \leq\eta_1\int\rho|\dot{\mathbf{u}}|^2dx
+C(\eta_1)\int|\mathbf{u}|^2|\nabla\mathbf{u}|^2dx
+\frac{d}{dt}\int\left[\left(P+\frac{|\mathbf{b}|^2}{2}\right)\divv\mathbf{u}
-(\mathbf{b}\otimes\mathbf{b}):\nabla\mathbf{u}\right]dx \nonumber\\
& \quad -\int\left(P+\frac{|\mathbf{b}|^2}{2}\right)_t\divv\mathbf{u}dx
+\int(\mathbf{b}\otimes\mathbf{b})_{t}:\nabla\mathbf{u}dx
\nonumber\\
& \leq \frac{d}{dt}\int\left[\left(P+\frac{|\mathbf{b}|^2}{2}\right)\divv\mathbf{u}
-(\mathbf{b}\otimes\mathbf{b}):\nabla\mathbf{u}
-\frac{\left(P+\frac{|\mathbf{b}|^2}{2}\right)^2}{2(\lambda+2\mu)}\right]dx
+\eta_1\int\rho|\dot{\mathbf{u}}|^2dx \nonumber\\
& \quad +C(\eta_1)\int|\mathbf{u}|^2|\nabla\mathbf{u}|^2dx
-\frac{1}{\lambda+2\mu}\int\left(P+\frac{|\mathbf{b}|^2}{2}\right)_tGdx
+\int(\mathbf{b}\otimes\mathbf{b})_{t}:\nabla\mathbf{u}dx.
\end{align}
It follows from \eqref{1.10} that $E$ satisfies
\begin{equation}\label{5.3}
\left(\rho E+\frac{|\mathbf{b}|^2}{2}\right)_t+\divv(\rho\mathbf{u}E)=\divv\mathbf{H}
\end{equation}
with
\begin{equation*}
\mathbf{H}=(\mathbf{u}\times\mathbf{b})\times\mathbf{b}
+\nu(\curl\mathbf{b})\times\mathbf{b}
+(2\mu\mathfrak{D}(\mathbf{u})+\lambda\divv\mathbf{u}\mathbb{I}_3)\mathbf{u}-P\mathbf{u}.
\end{equation*}
Then we infer from \eqref{3.13}, \eqref{5.3}, and \eqref{1.10}$_1$ that
\begin{align}\label{5.5}
-\int\left(P+\frac{|\mathbf{b}|^2}{2}\right)_tGdx
& =-\int \left(\rho E+\frac{|\mathbf{b}|^2}{2}\right)_tGdx+\frac12\int (\rho|\mathbf{u}|^2)_tG dx \nonumber\\
&=\int\divv(\rho \mathbf{u}E-\mathbf{H})Gdx+\frac12\int\rho_t|\mathbf{u}|^2Gdx
+\int\rho\mathbf{u}\cdot\mathbf{u}_tGdx \nonumber\\
&=-\int(\rho \mathbf{u}E-\mathbf{H})\cdot \nabla Gdx-\frac12\int\divv(\rho\mathbf{u})
|\mathbf{u}|^2Gdx+\int\rho\mathbf{u}\cdot\mathbf{u}_tGdx \nonumber\\
&=-\int(P\mathbf{u}-\mathbf{H})\cdot \nabla Gdx
+\frac12\int\rho\mathbf{u}\cdot\nabla(|\mathbf{u}|^2)Gdx
+\int\rho\mathbf{u}\cdot\mathbf{u}_tGdx \nonumber\\
& \triangleq\sum_{i=1}^{3}I_i.
\end{align}
From H{\"o}lder's inequality, Sobolev's inequality, \eqref{3.3}, and \eqref{3.5}, we have
\begin{align}\label{5.6}
I_1 & \leq \int (P|\mathbf{u}|+|\mathbf{H}|)|\nabla G|dx \nonumber\\
& \leq C\int (P|\mathbf{u}|+|\mathbf{u}||\mathbf{b}|^2
+|\mathbf{b}||\nabla\mathbf{b}|+|\mathbf{u}||\nabla\mathbf{u}|)|\nabla G|dx \nonumber\\
& \leq \eta_1\|\nabla G\|_{L^2}^2
+C(\eta_1)\left(\|P\|_{L^3}^2\|\mathbf{u}\|_{L^6}^2+\|\mathbf{u}\|_{L^6}^2\|\mathbf{b}\|_{L^6}^4
+\|\mathbf{b}\|_{L^3}^2\|\nabla\mathbf{b}\|_{L^6}^2
+\||\mathbf{u}||\nabla\mathbf{u}|\|_{L^2}^2\right) \nonumber\\
& \leq \eta_1\|\nabla G\|_{L^2}^2+C(\eta_1)\left(\|\nabla\mathbf{u}\|_{L^2}^2+\|\nabla\mathbf{b}\|_{H^1}^2
+\||\mathbf{u}||\nabla\mathbf{u}|\|_{L^2}^2\right).
\end{align}
By \eqref{3.14}, \eqref{3.2}, \eqref{3.3}, \eqref{3.5}, H{\"o}lder's inequality, and Sobolev's inequality, one gets
\begin{align}\label{5.7}
I_2 & \leq \int \rho|\mathbf{u}|^2|\nabla\mathbf{u}||G|dx\nonumber\\
& \leq C\int \rho|\mathbf{u}|^2|\nabla\mathbf{u}|(|\nabla\mathbf{u}|
+|\mathbf{b}|^2+P)dx\nonumber\\
& \leq C\int\left(|\mathbf{u}|^2|\nabla\mathbf{u}|^2+|\mathbf{u}|^2|\nabla\mathbf{u}||\mathbf{b}|^2
+\rho|\mathbf{u}|^2|\nabla\mathbf{u}|\right)dx \nonumber\\
& \leq C(\||\mathbf{u}||\nabla\mathbf{u}|\|_{L^2}^2
+\|\mathbf{u}\|_{L^6}\||\mathbf{u}||\nabla\mathbf{u}|\|_{L^2}\|\mathbf{b}\|_{L^6}^2
+\|\rho\|_{L^3}\|\mathbf{u}\|_{L^6}\||\mathbf{u}||\nabla\mathbf{u}|\|_{L^2})\nonumber\\
& \leq C\||\mathbf{u}||\nabla\mathbf{u}|\|_{L^2}^2+C\|\nabla\mathbf{u}\|_{L^2}^2.
\end{align}
Similarly to $I_2$, we find that
\begin{align}\label{5.8}
I_3 &= \int\left[\rho\mathbf{u}\cdot\dot{\mathbf{u}} -\rho\mathbf{u}\cdot(\mathbf{u}\cdot\nabla)\mathbf{u}\right] Gdx \nonumber\\
& \leq C\int \rho(|\mathbf{u}||\dot{\mathbf{u}}|+|\mathbf{u}|^2|\nabla\mathbf{u}|)
(|\nabla\mathbf{u}|+|\mathbf{b}|^2+P)dx  \nonumber\\
& \leq C(\|\sqrt{\rho}\dot{\mathbf{u}}\|_{L^2}\||\mathbf{u}||\nabla\mathbf{u}|\|_{L^2}
+\|\sqrt{\rho}\dot{\mathbf{u}}\|_{L^2}
\|\mathbf{u}\|_{L^6}\|\rho\|_{L^6}^2+\||\mathbf{u}||\nabla\mathbf{u}|\|_{L^2}^2
+\||\mathbf{u}||\nabla\mathbf{u}|\|_{L^2}
\|\mathbf{u}\|_{L^6}\|\rho\|_{L^3}) \nonumber\\
& \leq\eta_1\|\sqrt{\rho}\dot{\mathbf{u}}\|_{L^2}^2
+C(\eta_1)(\||\mathbf{u}||\nabla\mathbf{u}|\|_{L^2}^2
+\|\nabla\mathbf{u}\|_{L^2}^2).
\end{align}
Inserting \eqref{5.6}--\eqref{5.8} into \eqref{5.5}, we arrive at
\begin{align}\label{5.08}
-\int \left(P+\frac{|\mathbf{b}|^2}{2}\right)_tGdx
& \leq \eta_1\|\nabla G\|_{L^2}^2+\eta_1\|\sqrt{\rho}\dot{\mathbf{u}}\|_{L^2}^2 \nonumber \\
& \quad +C(\eta_1)
(\|\nabla\mathbf{u}\|_{L^2}^2+\|\nabla\mathbf{b}\|_{H^1}^2
+\||\mathbf{u}||\nabla\mathbf{u}|\|_{L^2}^2).
\end{align}
In view of \eqref{3.14}, we obtain that
\begin{equation*}
\Delta G=\divv(\rho\dot{\mathbf{u}}-\mathbf{b}\cdot\nabla\mathbf{b}).
\end{equation*}
Then from the standard elliptic estimates, \eqref{3.2}, and \eqref{3.5}, we deduce that
\begin{align}\label{5.9}
\|\nabla G\|_{L^2}^2
& \leq C(\|\rho\dot{\mathbf{u}}\|_{L^2}^2
+\|\mathbf{b}\cdot\nabla\mathbf{b}\|_{L^2}^2)  \nonumber\\
& \leq C(\|\sqrt{\rho}\dot{\mathbf{u}}\|_{L^2}^2
+\|\mathbf{b}\|_{L^3}^2\|\nabla\mathbf{b}\|_{L^6}^2)  \nonumber\\
& \leq C(\|\sqrt{\rho}\dot{\mathbf{u}}\|_{L^2}^2+\|\nabla\mathbf{b}\|_{H^1}^2),
\end{align}
which combined with \eqref{5.08} implies that
\begin{align}\label{5.10}
-\int \left(P+\frac{|\mathbf{b}|^2}{2}\right)_tGdx
& \leq C\eta_1\|\sqrt{\rho}\dot{\mathbf{u}}\|_{L^2}^2+C(\eta_1)
(\|\nabla\mathbf{u}\|_{L^2}^2+\|\nabla\mathbf{b}\|_{H^1}^2
+\||\mathbf{u}||\nabla\mathbf{u}|\|_{L^2}^2).
\end{align}

For the last term on the right-hand side of \eqref{5.2}, we obtain from H{\"o}lder's inequality and \eqref{3.5} that
\begin{align}\label{5.13}
\int(\mathbf{b}\otimes\mathbf{b})_t:\nabla\mathbf{u}dx
&\leq C\int|\mathbf{b}_t||\mathbf{b}||\nabla\mathbf{u}|dx \nonumber\\
&\leq C\|\mathbf{b}_t\|_{L^2}\|\mathbf{b}\|_{L^\infty}\|\nabla\mathbf{u}\|_{L^2} \nonumber\\
&\leq \tilde{\eta}\|\mathbf{b}_t\|_{L^2}^2+C(\tilde{\eta})\|\nabla\mathbf{u}\|_{L^2}^2.
\end{align}
Inserting \eqref{5.10} and \eqref{5.13} into \eqref{5.2} and choosing $\eta_1$ suitably small, we have
\begin{align}\label{5.14}
\frac{d}{dt}\int \Phi dx+\|\sqrt{\rho}\dot{\mathbf{u}}\|_{L^2}^2
& \leq
\tilde{\eta}\|\mathbf{b}_t\|_{L^2}^2+C_2\|\nabla^2\mathbf{b}\|_{L^2}^2 +C(\|\nabla\mathbf{u}\|_{L^2}^2+\|\nabla\mathbf{b}\|_{L^2}^2
+\||\mathbf{u}||\nabla\mathbf{u}|\|_{L^2}^2),
\end{align}
where
$$\Phi=\mu|\nabla \mathbf{u}|^2+(\mu+\lambda)(\divv\mathbf{u})^2-\left(2P+|\mathbf{b}|^2\right)\divv \mathbf{u}+2(\mathbf{b}\otimes\mathbf{b}):\nabla\mathbf{u}
+\frac{\left(P+\frac{|\mathbf{b}|^2}{2}\right)^2}{\lambda+2\mu}$$
satisfies
\begin{equation}\label{5.15}
\frac{\mu}{2}\|\nabla\mathbf{u}\|_{L^2}^2-C\leq \int\Phi dx\leq\mu\|\nabla\mathbf{u}\|_{L^2}^2+C
\end{equation}
due to \eqref{3.3} and \eqref{3.5}.

It follows from $\eqref{1.10}_4$, H{\"o}lder's inequality, and \eqref{3.5} that
\begin{align}\label{5.16}
 &\nu\frac{d}{dt}
\int|\nabla\mathbf{b}|^2dx+\int|\mathbf{b}_t|^2dx+\nu^2\int|\Delta\mathbf{b}|^2dx
 \nonumber\\
&= \int |\mathbf{b}_t-\nu\Delta\mathbf{b}|^2dx \nonumber\\
&=\int |\mathbf{b}\cdot\nabla\mathbf{u}-\mathbf{u}\cdot \nabla\mathbf{b}-\mathbf{b}\divv\mathbf{u}|^2dx \nonumber\\
&\leq C\int |\mathbf{b}|^2|\nabla\mathbf{u}|^2\mbox{d}x
+C\int |\mathbf{u}|^2|\nabla\mathbf{b}|^2dx \nonumber\\
&\leq C\|\mathbf{b}\|_{L^\infty}^2\|\nabla\mathbf{u}\|_{L^2}^2+C\|\mathbf{u}\|_{L^6}^2
\|\nabla\mathbf{b}\|_{L^3}^2\nonumber\\
&\leq C\|\nabla\mathbf{u}\|_{L^2}^2
+C\|\nabla\mathbf{u}\|_{L^2}^2\|\nabla\mathbf{b}\|_{L^2}\|\nabla\mathbf{b}\|_{H^1} \nonumber\\
&\leq C\|\nabla\mathbf{u}\|_{L^2}^2
+C\|\nabla\mathbf{u}\|_{L^2}^2\|\nabla\mathbf{b}\|_{L^2}
(\|\nabla\mathbf{b}\|_{L^2}+\|\nabla^2\mathbf{b}\|_{L^2}) \nonumber\\
&\leq \eta_2\|\nabla^2\mathbf{b}\|_{L^2}^2
+C(\eta_2)(\|\nabla\mathbf{u}\|_{L^2}^2
+\|\nabla\mathbf{b}\|_{L^2}^2)(\|\nabla\mathbf{u}\|_{L^2}^2+\|\nabla\mathbf{b}\|_{L^2}^2+1).
\end{align}
Noting that the standard $L^2$ estimate of elliptic system gives
$$\|\nabla^2\mathbf{b}\|_{L^2}\leq C_3\|\Delta\mathbf{b}\|_{L^2},$$
hence we deduce after choosing $\eta_2$ suitably small that
\begin{align}\label{5.17}
& 2\nu\frac{d}{dt}\int |\nabla\mathbf{b}|^2dx+2\|\mathbf{b}_t\|_{L^2}^2
+C_3^{-1}\nu^2\|\nabla^2\mathbf{b}\|_{L^2}^2 \nonumber\\
& \leq C(\|\nabla\mathbf{u}\|_{L^2}^2
+\|\nabla\mathbf{b}\|_{L^2}^2)(\|\nabla\mathbf{u}\|_{L^2}^2
+\|\nabla\mathbf{b}\|_{L^2}^2+1).
\end{align}
Then adding \eqref{5.17} to \eqref{5.14} and choosing $\tilde{\eta}$ small enough, we have
\begin{align*}
\frac{d}{dt}
&\int(\Phi+2C_4\nu|\nabla\mathbf{b}|^2)dx
+\frac12(\|\sqrt{\rho}\dot{\mathbf{u}}\|_{L^2}^2+\|\mathbf{b}_t\|_{L^2}^2
+\|\nabla^2\mathbf{b}\|_{L^2}^2)   \\
&\leq C(\|\nabla\mathbf{u}\|_{L^2}^2+\|\nabla\mathbf{b}\|_{L^2}^2)
(\|\nabla\mathbf{u}\|_{L^2}^2+\|\nabla\mathbf{b}\|_{L^2}^2+1)   \\
&~~~~+C(\|\nabla\mathbf{u}\|_{L^2}^2+\|\nabla\mathbf{b}\|_{L^2}^2
+\||\mathbf{u}||\nabla\mathbf{u}|\|_{L^2}^2).
\end{align*}
Then we obtain the desired \eqref{5.1} after using Gronwall's inequality, \eqref{3.3}, \eqref{4.1}, and \eqref{5.15}. This completes the proof of Lemma \ref{lem35}.            \hfill $\Box$

Next, we have the following estimates on the material derivatives of the velocity which are important for the higher order estimates of strong solutions.
\begin{lemma}\label{lem36}
Under the condition \eqref{3.1}, it holds that for any $T\in[0,T^*)$,
\begin{equation}\label{6.1}
\sup_{0\leq t\leq T}\left(\|\sqrt{\rho}\dot{\mathbf{u}}\|_{L^2}^2
+\|\mathbf{b}_t\|_{L^2}^2+\|\nabla^2\mathbf{b}\|_{L^2}^2\right)
+\int_0^T\left(\|\nabla\dot{\mathbf{u}}\|_{L^2}^2+\|\nabla\mathbf{b}_t\|_{L^2}^2\right)dt \leq C.
\end{equation}
\end{lemma}
{\it Proof.}
By the definition of $\dot{\mathbf{u}}$, we can rewrite $\eqref{1.10}_2$ as follows:
\begin{equation}\label{6.2}
\rho\dot{\mathbf{u}}+\nabla P=\mu\Delta\mathbf{u}
+(\lambda+\mu)\nabla\divv\mathbf{u}+\curl\mathbf{b}\times\mathbf{b}.
\end{equation}
Differentiating \eqref{6.2} with respect to $t$ and using \eqref{1.10}$_1$, we have
\begin{align}\label{6.3}
\rho\dot{\mathbf{u}}_t+\rho\mathbf{u}\cdot\nabla\dot{\mathbf{u}}
+\nabla P_t&=\mu\Delta\dot{\mathbf{u}}+(\lambda+\mu)\divv\dot{\mathbf{u}}
-\mu\Delta(\mathbf{u}\cdot\nabla\mathbf{u})
-(\lambda+\mu)\divv(\mathbf{u}\cdot\nabla\mathbf{u})\nonumber\\
& \quad +(\curl\mathbf{b}\times\mathbf{b})_t+\divv(\rho\dot{\mathbf{u}}\otimes\mathbf{u}).
\end{align}
Multiplying \eqref{6.3} by $\dot{\mathbf{u}}$ and integrating by parts over $\mathbb{R}^3$, we get
\begin{align}\label{6.4}
&\frac12\frac{d}{dt}\int \rho|\dot{\mathbf{u}}|^2\mbox{d}x+\mu\int|\nabla\dot{\mathbf{u}}|^2dx
+(\lambda+\mu)\int|\divv\dot{\mathbf{u}}|^2dx  \nonumber\\
&=\int P_t\divv\dot{\mathbf{u}}+(\nabla P\otimes\mathbf{u}):\nabla\dot{\mathbf{u}}dx
-\int [\divv(\curl\mathbf{b}\times\mathbf{b})\otimes\mathbf{u}
-(\curl\mathbf{b}\times\mathbf{b})_t]\cdot \dot{\mathbf{u}}dx   \nonumber\\
& \quad +\mu\int [\divv(\Delta\mathbf{u}\otimes\mathbf{u})-\Delta(\mathbf{u}\cdot\nabla\mathbf{u})]\cdot \dot{\mathbf{u}}dx+(\lambda+\mu)\int [(\nabla\divv\mathbf{u})\otimes\mathbf{u}
-\nabla\divv(\mathbf{u}\cdot\nabla\mathbf{u})]\cdot \dot{\mathbf{u}}dx   \nonumber\\
&\triangleq\sum_{i=1}^{4}J_i,
\end{align}
where $J_i$ can be bounded as follows.

It follows from $\eqref{1.10}_3$ that
\begin{align}\label{6.5}
J_1 & =\int \left(-\divv(P\mathbf{u})\divv\dot{\mathbf{u}}
-P\divv\mathbf{u}\divv\dot{\mathbf{u}}
+\mathcal{T}(\mathbf{u}):\nabla\mathbf{u}\divv\dot{\mathbf{u}}
+\nu|\curl\mathbf{b}|^2\divv\dot{\mathbf{u}}+(\nabla P\otimes\mathbf{u}):\nabla\dot{\mathbf{u}}\right)dx   \nonumber\\
&=\int \left(P\mathbf{u}\nabla\divv\dot{\mathbf{u}}-P\divv\mathbf{u}\divv\dot{\mathbf{u}}
+\mathcal{T}(\mathbf{u}):\nabla\mathbf{u}\divv\dot{\mathbf{u}}
+\nu|\curl\mathbf{b}|^2\divv\dot{\mathbf{u}}\right)dx \nonumber\\
& \quad -\int\left(P\nabla\mathbf{u}^\top:\nabla\dot{\mathbf{u}}
+P\mathbf{u}\nabla\divv\dot{\mathbf{u}}\right)dx  \nonumber\\
&=\int \left(-P\divv\mathbf{u}\divv\dot{\mathbf{u}}
+\mathcal{T}(\mathbf{u}):\nabla\mathbf{u}\divv\dot{\mathbf{u}}
+\nu|\curl\mathbf{b}|^2\divv\dot{\mathbf{u}}
-P\nabla\mathbf{u}^\top:\nabla\dot{\mathbf{u}}\right)dx\nonumber\\
&\leq C\int \left(|\nabla\mathbf{u}||\nabla\dot{\mathbf{u}}|
+|\nabla\mathbf{u}|^2|\nabla\dot{\mathbf{u}}|
+|\nabla\mathbf{b}|^2|\nabla\dot{\mathbf{u}}|\right)dx  \nonumber\\
&\leq C(\|\nabla\mathbf{u}\|_{L^2}+\|\nabla\mathbf{u}\|_{L^4}^2
+\|\nabla\mathbf{b}\|_{L^4}^2)\|\nabla\dot{\mathbf{u}}\|_{L^2}  \nonumber\\
&\leq C(\|\nabla\mathbf{u}\|_{L^2}+\|\nabla\mathbf{u}\|_{L^4}^2
+\|\nabla\mathbf{b}\|_{L^4}^2)\|\nabla\dot{\mathbf{u}}\|_{L^2},
\end{align}
where $\mathcal{T}(\mathbf{u})=2\mu\mathfrak{D}(\mathbf{u})+\lambda\divv\mathbf{u}\mathbb{I}_3$.
Integrating by parts leads to
\begin{align}\label{6.6}
J_2&=\int \left[\divv(\mathbf{b}\otimes\mathbf{b})_t-\nabla\left(\frac{|\mathbf{b}|^2}{2}\right)_t
-\divv(\mbox{curl}\mathbf{b}\times\mathbf{b})\otimes\mathbf{u}\right]\cdot \dot{\mathbf{u}}\mbox{d}x\nonumber\\
& \leq C\int(|\mathbf{b}||\mathbf{b}_t|+|\mathbf{b}||\nabla\mathbf{b}||\mathbf{u}|)
|\nabla\dot{\mathbf{u}}|\mbox{d}x\nonumber\\
&\leq C(\|\mathbf{b}\|_{L^6}\|\mathbf{b}_t\|_{L^3}
+\|\mathbf{b}\|_{L^6}\|\nabla\mathbf{b}\|_{L^6}\|\mathbf{u}\|_{L^6})
\|\nabla\dot{\mathbf{u}}\|_{L^2}\nonumber\\
&\leq C(\|\mathbf{b}_t\|_{L^2}^{\frac12}\|\nabla\mathbf{b}_t\|_{L^2}^{\frac12}
+\|\nabla\mathbf{b}\|_{H^1})\|\nabla\dot{\mathbf{u}}\|_{L^2}.
\end{align}
For $J_3$ and $J_4$, notice that for all $1\leq i,j,k\leq 3,$ one has
\begin{align}
\partial_j(\partial_{kk}u_iu_j)-\partial_{kk}(u_j\partial_j u_i)&=\partial_k(\partial_ju_j\partial_ku_i)-\partial_k(\partial_ku_j\partial_ju_i)
-\partial_j(\partial_ku_j\partial_ku_i),\nonumber\\
\partial_j(\partial_{ik}u_ku_j)-\partial_{ij}(u_k\partial_ku_j)
&=\partial_i(\partial_ju_j\partial_ku_k)-\partial_i(\partial_ju_k\partial_ku_j)
-\partial_k(\partial_iu_k\partial_ju_j).\nonumber
\end{align}
So integrating by parts gives
\begin{align}
J_3&=\mu\int [\partial_k(\partial_ju_j\partial_ku_i)-\partial_k(\partial_ku_j\partial_ju_i)
-\partial_j(\partial_ku_j\partial_ku_i)]\dot{u_i}\mbox{d}x\nonumber\\
&\leq C\|\nabla\mathbf{u}\|_{L^4}^2\|\nabla\dot{\mathbf{u}}\|_{L^2},\label{6.7}\\
J_4&=(\lambda+\mu)\int [\partial_i(\partial_ju_j\partial_ku_k)-\partial_i(\partial_ju_k\partial_ku_j)
-\partial_k(\partial_iu_k\partial_ju_j)]\dot{u_i}\mbox{d}x\nonumber\\
&\leq C\|\nabla\mathbf{u}\|_{L^4}^2\|\nabla\dot{\mathbf{u}}\|_{L^2}.\label{6.8}
\end{align}
Inserting \eqref{6.5}--\eqref{6.8} into \eqref{6.4} and applying \eqref{5.1} lead to
\begin{align}\label{6.9}
& \frac12\frac{d}{dt}
\|\sqrt{\rho}\dot{\mathbf{u}}\|_{L^2}^2+\mu\|\nabla\dot{\mathbf{u}}\|_{L^2}^2
+(\lambda+\mu)\|\divv\dot{\mathbf{u}}\|_{L^2}^2  \nonumber\\
& \leq C(\|\nabla\mathbf{u}\|_{L^2}+\|\nabla\mathbf{u}\|_{L^4}^2
+\|\nabla\mathbf{b}\|_{L^4}^2+\|\mathbf{b}_t\|_{L^2}^{\frac12}\|\nabla\mathbf{b}_t\|_{L^2}^{\frac12}
+\|\nabla\mathbf{b}\|_{H^1})\|\nabla\dot{\mathbf{u}}\|_{L^2}\nonumber\\
&\leq \delta_1\|\nabla\dot{\mathbf{u}}\|_{L^2}^2+\delta_2\|\nabla\mathbf{b}_t\|_{L^2}^2 +C(\delta_1,\delta_2)(\|\nabla\mathbf{u}\|_{L^4}^4+\|\nabla\mathbf{b}\|_{L^4}^4
+\|\mathbf{b}_t\|_{L^2}^2+\|\nabla^2\mathbf{b}\|_{L^2}^2+1).
\end{align}

From $\eqref{1.10}_4$, the standard regularity estimate of elliptic equations to \eqref{1.10}$_4$, \eqref{3.5}, and \eqref{5.1}, we get
\begin{align*}
\|\nabla^2\mathbf{b}\|_{L^2}^2&\leq C(\|\mathbf{b}_t\|_{L^2}^2
+\||\mathbf{u}||\nabla\mathbf{b}|\|_{L^2}^2+\||\mathbf{b}||\nabla\mathbf{u}|\|_{L^2}^2)\nonumber\\
&\leq C(\|\mathbf{b}_t\|_{L^2}^2+\|\mathbf{u}\|_{L^6}^2\|\nabla\mathbf{b}\|_{L^3}^2
+\|\mathbf{b}\|_{L^{\infty}}^2\|\nabla\mathbf{u}\|_{L^2}^2)\nonumber\\
&\leq C(\|\mathbf{b}_t\|_{L^2}^2+\|\nabla\mathbf{u}\|_{L^2}^2
\|\nabla\mathbf{b}\|_{L^2}\|\nabla\mathbf{b}\|_{H^1}
+1) \nonumber\\
&\leq \frac{1}{2}\|\nabla^2\mathbf{b}\|_{L^2}^2+C\|\mathbf{b}_t\|_{L^2}^2+C,
\end{align*}
which implies that
\begin{equation}\label{6.10}
\|\nabla^2\mathbf{b}\|_{L^2}^2\leq C\|\mathbf{b}_t\|_{L^2}^2+C.
\end{equation}
Differentiating $\eqref{1.10}_4$ with respect to $t$, we have
\begin{equation}\label{6.11}
\mathbf{b}_{tt}-\nu\Delta\mathbf{b}_t
=\mathbf{b}_t\cdot\nabla\mathbf{u}
-\mathbf{u}\cdot\nabla\mathbf{b}_t
-\mathbf{b}_t\divv\mathbf{u}+\mathbf{b}\cdot\nabla\mathbf{u}_t
-\mathbf{u}_t\cdot\nabla\mathbf{b}-\mathbf{b}\divv\mathbf{u}_t.
\end{equation}
Multiplying \eqref{6.11} by $\mathbf{b}_t$ and integrating by parts lead to
\begin{align}\label{6.12}
\frac12\frac{d}{dt}\int|\mathbf{b}_t|^2dx+\nu\int|\nabla\mathbf{b}_t|^2dx
&=\int(\mathbf{b}_t\cdot\nabla\mathbf{u}-\mathbf{u}\cdot\nabla\mathbf{b}_t
-\mathbf{b}_t\divv\mathbf{u})\cdot\mathbf{b}_tdx \nonumber\\
& \quad +\int(\mathbf{b}\cdot\nabla\mathbf{u}_t-\mathbf{u}_t\cdot\nabla\mathbf{b}
-\mathbf{b}\divv\mathbf{u}_t)\cdot\mathbf{b}_tdx\nonumber\\
&\triangleq L_1+L_2.
\end{align}
Integrating by parts implies that
\begin{align}\label{6.13}
L_1&=\int \left(\mathbf{b}_t\cdot\nabla\mathbf{u}\cdot\mathbf{b}_t
-\frac12|\mathbf{b}_t|^2\divv\mathbf{u}\right)dx  \nonumber\\
&\leq C\|\mathbf{b}_t\|_{L^4}^2\|\nabla\mathbf{u}\|_{L^2}  \nonumber\\
&\leq C\|\mathbf{b}_t\|_{L^2}^{\frac12}\|\nabla\mathbf{b}_t\|_{L^2}^{\frac32} \nonumber\\
&\leq \delta_1\|\nabla\mathbf{b}_t\|_{L^2}^2+C(\delta_1)\|\mathbf{b}_t\|_{L^2}^2,
\end{align}
and
\begin{align}\label{6.14}
L_2 & =\int (\mathbf{b}\cdot\nabla\dot{\mathbf{u}}-\dot{\mathbf{u}}\cdot\nabla\mathbf{b}
-\mathbf{b}\divv\dot{\mathbf{u}})\cdot\mathbf{b}_tdx   \nonumber\\
& \quad -\int (\mathbf{b}\cdot\nabla(\mathbf{u}\cdot\nabla\mathbf{u})
-(\mathbf{u}\cdot\nabla\mathbf{u})\cdot\nabla\mathbf{b}
-\mathbf{b}\divv(\mathbf{u}\cdot\nabla\mathbf{u}))\cdot\mathbf{b}_t   \nonumber\\
&= \int (\mathbf{b}\cdot\nabla\dot{\mathbf{u}}-\dot{\mathbf{u}}\cdot\nabla\mathbf{b}
-\mathbf{b}\divv\dot{\mathbf{u}})\cdot\mathbf{b}_t\mbox{d}x   \nonumber\\
& \quad +\int (\mathbf{u}\cdot\nabla\mathbf{u})\cdot(\mathbf{b}\cdot\nabla\mathbf{b}_t)
+(\mathbf{u}\cdot\nabla\mathbf{u})\cdot\nabla\mathbf{b}_t\cdot\mathbf{b}dx   \nonumber\\
&\leq C\int |\mathbf{b}||\mathbf{b}_t||\nabla\dot{\mathbf{u}}|+|\dot{\mathbf{u}}||\nabla\mathbf{b}||\mathbf{b}_t|
+|\mathbf{u}||\nabla\mathbf{u}||\mathbf{b}||\nabla\mathbf{b}_t|dx  \nonumber\\
&\leq C(\|\mathbf{b}\|_{L^6}\|\mathbf{b}_t\|_{L^3}\|\nabla\dot{\mathbf{u}}\|_{L^2}
+\|\dot{\mathbf{u}}\|_{L^6}\|\nabla\mathbf{b}\|_{L^2}\|\mathbf{b}_t\|_{L^3}
+\|\mathbf{u}\|_{L^6}\|\nabla\mathbf{u}\|_{L^6}\|\mathbf{b}\|_{L^6}
\|\nabla\mathbf{b}_t\|_{L^2})  \nonumber\\
&\leq C(\|\mathbf{b}_t\|_{L^3}\|\nabla\dot{\mathbf{u}}\|_{L^2}
+\|\nabla\mathbf{u}\|_{L^6}\|\nabla\mathbf{b}_t\|_{L^2})   \nonumber\\
&\leq C(\|\mathbf{b}_t\|_{L^2}^{\frac12}\|\nabla\mathbf{b}_t\|_{L^2}^{\frac12}
\|\nabla\dot{\mathbf{u}}\|_{L^2}
+\|\nabla\mathbf{u}\|_{L^6}\|\nabla\mathbf{b}_t\|_{L^2})   \nonumber\\
&\leq \delta_1\|\nabla\mathbf{b}_t\|_{L^2}^2
+\delta_2\|\nabla\dot{\mathbf{u}}\|_{L^2}^2
+C(\delta_1,\delta_2)\|\mathbf{b}_t\|_{L^2}+C(\delta_1)\|\nabla\mathbf{u}\|_{L^6}^2.
\end{align}
Inserting \eqref{6.13} and \eqref{6.14} into \eqref{6.12}, we have
\begin{align}\label{6.15}
\frac12\frac{d}{dt}\|\mathbf{b}_t\|_{L^2}^2+\nu\|\nabla\mathbf{b}_t\|_{L^2}^2
\leq 2\delta_1\|\nabla\mathbf{b}_t\|_{L^2}^2
+\delta_2\|\nabla\dot{\mathbf{u}}\|_{L^2}^2
+C(\delta_1,\delta_2)\|\mathbf{b}_t\|_{L^2}+C(\delta_1)\|\nabla\mathbf{u}\|_{L^6}^2.
\end{align}

Adding \eqref{6.15} to \eqref{6.9} and applying \eqref{6.10}, we obtain after choosing $\delta_1,\delta_2$ suitably small that
\begin{align}\label{6.16}
&\frac{d}{dt}\left(\|\sqrt{\rho}\dot{\mathbf{u}}\|_{L^2}^2
+\|\mathbf{b}_t\|_{L^2}^2\right)
+\tilde{C}\left(\|\nabla\dot{\mathbf{u}}\|_{L^2}^2+\|\nabla\mathbf{b}_t\|_{L^2}^2\right)
\nonumber\\
&\leq C(\|\mathbf{b}_t\|_{L^2}^2+\|\nabla\mathbf{u}\|_{L^4}^4
+\|\nabla\mathbf{b}\|_{L^4}^4+\|\nabla\mathbf{u}\|_{L^6}^2+1).
\end{align}
To estimate $\|\nabla\mathbf{u}\|_{L^6},$ let $\mathbf{u}=\mathbf{v}+\mathbf{w}$ such that
\begin{equation*}
\begin{cases}
\mu\Delta\mathbf{v}+(\lambda+\mu)\nabla\divv\mathbf{v}
=\nabla\left(P+\frac{|\mathbf{b}|^2}{2}\right),\\
\mathbf{v}(x,t)\to 0,~~~\text{as}~|x|\to +\infty;
\end{cases}
\end{equation*}
and
\begin{equation*}
\begin{cases}
\mu\Delta\mathbf{w}+(\lambda+\mu)\nabla\divv\mathbf{w}=\rho\dot{\mathbf{u}}
-\mathbf{b}\cdot\nabla\mathbf{b},\\
\mathbf{w}(x,t)\to 0,~~~\text{as}~|x|\to +\infty,
\end{cases}
\end{equation*}
which implies that
\begin{equation*}
\|\nabla\mathbf{v}\|_{L^6}\leq C\left\|P+\frac{|\mathbf{b}|^2}{2}\right\|_{L^6}\leq C,
\end{equation*}
and
\begin{align*}
\|\nabla\mathbf{w}\|_{L^6}+\|\nabla^2\mathbf{w}\|_{L^2}
\leq C(\|\rho\dot{\mathbf{u}}\|_{L^2}+\|\mathbf{b}\cdot\nabla\mathbf{b}\|_{L^2})
\leq C(\|\sqrt{\rho}\dot{\mathbf{u}}\|_{L^2}+\|\nabla\mathbf{b}\|_{H^1}).
\end{align*}
Then we have
\begin{align}\label{5.12}
\|\nabla\mathbf{u}\|_{L^6}^2
\leq\|\nabla\mathbf{v}\|_{L^6}^2+\|\nabla\mathbf{w}\|_{L^6}^2
\leq C(\|\sqrt{\rho}\dot{\mathbf{u}}\|_{L^2}^2+\|\nabla\mathbf{b}\|_{H^1}^2)+C
\end{align}
By \eqref{5.12}, \eqref{5.1}, and \eqref{6.10}, one has
\begin{align}\label{6.17}
\|\nabla\mathbf{u}\|_{L^6}^2&\leq C(\|\sqrt{\rho}\dot{\mathbf{u}}\|_{L^2}^2+\|\nabla\mathbf{b}\|_{H^1}^2+1)\leq C(\|\sqrt{\rho}\dot{\mathbf{u}}\|_{L^2}^2+\|\mathbf{b}_t\|_{L^2}^2+1).
\end{align}
It follows from H{\"o}lder's inequality, \eqref{5.1}, and \eqref{6.17} that
\begin{align}\label{6.18}
\|\nabla\mathbf{u}\|_{L^4}^4\leq C\|\nabla\mathbf{u}\|_{L^2}\|\nabla\mathbf{u}\|_{L^6}^3
& \leq C(\|\sqrt{\rho}\dot{\mathbf{u}}\|_{L^2}^3+\|\mathbf{b}_t\|_{L^2}^3+1)
\nonumber \\
& \leq C\left(1+\|\sqrt{\rho}\dot{\mathbf{u}}\|_{L^2}^2+\|\mathbf{b}_t\|_{L^2}^2\right)
\left(\|\sqrt{\rho}\dot{\mathbf{u}}\|_{L^2}^2+\|\mathbf{b}_t\|_{L^2}^2\right)+C.
\end{align}
Similarly, we get
\begin{align}\label{6.19}
\|\nabla\mathbf{b}\|_{L^4}^4
& \leq C\|\nabla\mathbf{b}\|_{L^2}\|\nabla\mathbf{b}\|_{H^1}^3 \nonumber \\
& \leq C(\|\nabla^2\mathbf{b}\|_{L^2}^3+1) \nonumber \\
& \leq C(\|\mathbf{b}_t\|_{L^2}^3+1) \nonumber \\
& \leq C(1+\|\mathbf{b}_t\|_{L^2}^2)\|\mathbf{b}_t\|_{L^2}^2+C.
\end{align}
Substituting \eqref{6.17}--\eqref{6.19} into \eqref{6.16} and then
applying Gronwall's inequality and \eqref{5.1} give the desired \eqref{6.1}. Hence we complete the proof of Lemma \ref{lem36}.
\hfill $\Box$

Finally, the following lemma will treat the higher order derivatives of the solutions which are needed to guarantee the extension of local strong solution to be a global one.
\begin{lemma}\label{lem37}
Under the condition \eqref{3.1}, and let $\tilde{q}\in(3,6]$ be as in Theorem \ref{thm1.1}, then it holds that for any $T\in[0,T^*)$,
\begin{equation}\label{7.1}
\sup_{0\leq t\leq T}\left(\|(\rho,P)\|_{H^1\cap W^{1,\tilde{q}}}+\|\nabla\mathbf{u}\|_{H^1}
+\|\mathbf{b}\|_{H^2}\right) \leq C.
\end{equation}
\end{lemma}
{\it Proof.}
First, in view of \eqref{3.3}, \eqref{5.1}, and \eqref{6.1}, one has
\begin{align}\label{7.2}
\|\mathbf{b}\|_{H^2} \leq C.
\end{align}
It follows from \eqref{5.12}, \eqref{6.1}, and \eqref{7.2} that
\begin{align}\label{7.3}
\|\nabla\mathbf{u}\|_{L^6} \leq C.
\end{align}
By virtue of Gagliardo-Nirenberg inequality, Sobolev's inequality, \eqref{5.1}, and \eqref{7.3}, we arrive at
\begin{align}\label{7.4}
\|\mathbf{u}\|_{L^{\infty}}
\leq C\|\mathbf{u}\|_{L^6}^{\frac12}
\|\nabla\mathbf{u}\|_{L^6}^{\frac12}
\leq C\|\nabla\mathbf{u}\|_{L^2}^{\frac12}
\|\nabla\mathbf{u}\|_{L^6}^{\frac12}\leq C.
\end{align}
For $2\leq p\leq\tilde{q}\leq6$, direct calculations show that
\begin{equation}\label{7.5}
\frac{d}{dt}\|\nabla\rho\|_{L^p}
\leq C(1+\|\nabla\mathbf{u}\|_{L^{\infty}})\|\nabla\rho\|_{L^p}
+C\|\nabla^2\mathbf{u}\|_{L^p}.
\end{equation}
Similarly,
\begin{equation}\label{7.6}
\frac{d}{dt}\|\nabla P\|_{L^p}\leq C(1+\|\nabla\mathbf{u}\|_{L^{\infty}})(\|\nabla P\|_{L^p}+\|\nabla^2\mathbf{u}\|_{L^p})+C\|\nabla\mathbf{b}\|_{L^{\infty}}
\|\nabla^2\mathbf{b}\|_{L^p}.
\end{equation}
Applying the standard $L^p$-estimate of elliptic system to \eqref{3.15}, \eqref{3.2}, and \eqref{7.2} yield
\begin{equation*}
\|\nabla G\|_{L^6}+\|\nabla\pmb\omega\|_{L^6}
\leq C(\|\rho\dot {\mathbf{u}}\|_{L^6}+\|\mathbf{b}\cdot\nabla\mathbf{b}\|_{L^6})
\leq C+C\|\nabla\dot {\mathbf{u}}\|_{L^2} ,
\end{equation*}
which combined with Gagliardo-Nirenberg inequality implies
\begin{align*}
\|G\|_{L^\infty}
&\leq\|G\|_{L^2}^{\beta}\|\nabla G\|_{L^6}^{1-\beta}
\leq C+ C\|\nabla\dot {\mathbf{u}}\|_{L^2}^{1-\beta},\\
\|\pmb\omega\|_{L^\infty}
&\leq\|\pmb\omega\|_{L^2}^{\beta}\|\nabla \pmb\omega\|_{L^6}^{1-\beta}
\leq C+ C\|\nabla\dot {\mathbf{u}}\|_{L^2}^{1-\beta},
\end{align*}
for some $\beta\in (0,1)$.

Employing the standard $L^p$-estimate of elliptic system to \eqref{1.10}$_2$ leads to
\begin{align}\label{7.7}
\|\nabla^2\mathbf{u}\|_{L^p}
& \leq C\left(\|\rho\dot{\mathbf{u}}\|_{L^p}+\|\nabla P\|_{L^p}
+\|\mathbf{b}\cdot\nabla\mathbf{b}\|_{L^p}\right) \nonumber\\
&\leq C\left(1+\|\rho\dot{\mathbf{u}}\|_2^\alpha\|\rho\dot {\mathbf{u}}\|_{L^6}^{1-\alpha}+\|\nabla P\|_{L^p}\right) \nonumber\\
&\leq C\left(1+\|\rho\dot{\mathbf{u}}\|_2^\alpha\|\nabla\dot {\mathbf{u}}\|_{L^2}^{1-\alpha}+\|\nabla P\|_{L^p}\right) \nonumber\\
&\leq C\left(1+\|\nabla\dot{\mathbf{u}}\|_{L^2}^{1-\alpha}+\|\nabla P\|_{L^p}\right),
\end{align}
for some $\alpha\in (0,1)$.
This together with Lemma \ref{lem23} gives
\begin{equation}\label{7.8}
\|\nabla\mathbf{u}\|_{L^\infty}
\leq C\left(1+\|\nabla\dot{\mathbf{u}}\|_{L^2}^{1-\beta}\right)
\log \left(e+\|\nabla\dot{\mathbf{u}}\|_{L^2} + \|\nabla P\|_{L^{p}}\right)
+ C\|\nabla\dot {\mathbf{u}}\|_{L^2}.
\end{equation}
Applying the standard $L^p$-estimate to \eqref{1.10}$_4$ yields
\begin{align}\label{7.9}
\|\nabla^2 \mathbf{b}\|_{L^p}&\leq C\left(\|\mathbf{b}_t\|_{L^p}
+\||\mathbf{u}||\nabla\mathbf{b}|\|_{L^p}+\||\mathbf{b}||\nabla\mathbf{u}|\|_{L^p}\right)
\nonumber\\
&\leq C\left(\|\mathbf{b}_t\|_{L^2}^{\frac{6-p}{2p}}
\|\nabla\mathbf{b}_t\|_{L^2}^{\frac{3p-6}{2p}}
+\|\mathbf{u}\|_{L^\infty}\|\nabla\mathbf{b}\|_{L^p}
+\|\mathbf{b}\|_{L^\infty}\|\nabla\mathbf{u}\|_{L^p}\right)\nonumber\\
&\leq C\left(\|\mathbf{b}_t\|_{L^2}^{\frac{6-p}{2p}}\|\nabla\mathbf{b}_t\|_{L^2}^{\frac{3p-6}{2p}}
+\|\nabla\mathbf{b}\|_{L^2}^{\frac{6-p}{2p}}
\|\nabla^2\mathbf{b}\|_{L^2}^{\frac{3p-6}{2p}}+\|\nabla\mathbf{u}\|_{L^2}^{\frac{6-p}{2p}}
\|\nabla\mathbf{u}\|_{L^6}^{\frac{3p-6}{2p}}\right)\nonumber\\
&\leq C\left(1+\|\nabla\mathbf{b}_t\|_{L^2}^{\frac{3p-6}{2p}}\right).
\end{align}
It follows from Lemma \ref{lem22} that
\begin{equation}\label{7.10}
\|\nabla\mathbf{b}\|_{L^\infty}\leq C(\|\nabla^2\mathbf{b}\|_{L^q}+1).
\end{equation}
Substituting \eqref{7.9} and \eqref{7.10} into \eqref{7.5}--\eqref{7.6} yields that
\begin{equation}\label{7.11}
f'(t)\leq  Cg(t)f(t)\log{f(t)}+Cg(t)f(t)+Cg(t),
\end{equation}
where
\begin{align*}
f(t)&\triangleq e+\|\nabla\rho\|_{L^{ q}}+\|\nabla P\|_{L^{ q}},\\
g(t)&\triangleq (1+\|\nabla\dot{\mathbf{u}}\|_{L^2})
\log(e+\|\nabla\dot{\mathbf{u}}\|_{L^2})+\|\nabla\mathbf{b}_t\|_{L^2}^2.
\end{align*}
This yields
\begin{equation}\label{7.12}
(\log f(t))'\leq Cg(t)+Cg(t)\log f(t)
\end{equation}
due to $f(t)>1.$
Thus it follows from  \eqref{7.12}, \eqref{6.1}, and  Gronwall's inequality that
\begin{equation}\label{7.13}
\sup_{0\leq t\leq T}\|(\nabla\rho,\nabla P)\|_{L^{q}}\leq C,
\end{equation}
which, combined with \eqref{7.8} and \eqref{6.1} gives that
\begin{equation}\label{7.14}
\int_{0}^{T}\|\nabla\mathbf{u}\|^2_{L^\infty}dt\leq C.
\end{equation}
Taking $p=2$ in \eqref{7.5}, one can get by using \eqref{7.14}, \eqref{7.7} and Gronwall's inequality  that
\begin{equation}\label{7.15}
\sup_{0\leq t\leq T}\|(\nabla\rho,\nabla P)\|_{L^2}\leq C,
\end{equation}
which together with \eqref{7.7} yields that
\begin{equation*}
\sup_{0\leq t\leq T}\|\nabla^2\mathbf{u}\|_{L^2}\leq C.
\end{equation*}
This combined  with \eqref{7.13}, \eqref{7.15}, \eqref{5.1}, and \eqref{7.2} finishes the proof of Lemma \ref{lem37}.
\hfill $\Box$

With Lemmas \ref{lem31}--\ref{lem37} at hand, we are now in a position to prove Theorem \ref{thm1.1}.

\textbf{Proof of Theorem \ref{thm1.1}.}
We argue by contradiction. Suppose that \eqref{B} were false, that is, \eqref{3.1} holds. Note that the general constant $C$ in Lemmas \ref{lem31}--\ref{lem37} is independent of $t<T^{*}$, that is, all the a priori estimates obtained in Lemmas \ref{lem31}--\ref{lem37} are uniformly bounded for any $t<T^{*}$. Hence, the function
\begin{equation*}
(\rho,\mathbf{u},P,\mathbf{b})(x,T^{*})
\triangleq\lim_{t\rightarrow T^{*}}(\rho,\mathbf{u},P,\mathbf{b})(x,t)
\end{equation*}
satisfy the initial condition \eqref{A} at $t=T^{*}$.

Furthermore, standard arguments yield that $\rho\dot{\mathbf{u}}\in C([0,T];L^2)$, which
implies $$ \rho\dot{\mathbf{u}}(x,T^\ast)=\lim_{t\rightarrow
T^\ast}\rho\dot{\mathbf{u}}\in L^2. $$
Hence, 
$$-\mu\Delta{\mathbf{u}}-(\lambda+\mu)\nabla\mbox{div}\mathbf{u}+\nabla P-\mbox{curl}\mathbf{b}\times\mathbf{b}
|_{t=T^\ast}=\sqrt{\rho}(x,T^\ast)g(x)
$$ with 
$$g(x)\triangleq
\begin{cases}
\rho^{-1/2}(x,T^\ast)(\rho\dot{\mathbf{u}})(x,T^\ast),&
\mbox{for}~~x\in\{x|\rho(x,T^\ast)>0\},\\
0,&\mbox{for}~~x\in\{x|\rho(x,T^\ast)=0\},
\end{cases}
$$
satisfying $g\in L^2$ due to \eqref{7.1}.
Therefore, one can take $(\rho,\mathbf{u},P,\mathbf{b})(x,T^\ast)$ as
the initial data and extend the local
strong solution beyond $T^\ast$. This contradicts the assumption on
$T^{\ast}$.

Thus we finish the proof of Theorem \ref{thm1.1}.
\hfill $\Box$

\section*{Acknowledgments}
The author would like to express his gratitude to the reviewers for careful reading and helpful suggestions which led to an improvement of the original manuscript.

\end{document}